\documentclass[11pt]{article}
\usepackage[a4paper,margin=2.6cm]{geometry}
\usepackage{amsmath,amssymb,amsthm,mathtools}
\usepackage{graphicx,booktabs,xcolor}
\usepackage{tikz,pgfplots}
\pgfplotsset{compat=1.16}
\usepackage{cite}
\usepackage[colorlinks=true,linkcolor=blue!50!black,citecolor=blue!50!black]{hyperref}

\newtheorem{theorem}{Theorem}[section]
\newtheorem{proposition}[theorem]{Proposition}
\newtheorem{lemma}[theorem]{Lemma}

\theoremstyle{definition}
\newtheorem{remark}[theorem]{Remark}

\renewcommand{\P}{\mathbb{P}}
\newcommand{\E}{\mathbb{E}}
\newcommand{\dd}{\mathrm{d}}

\newcommand{\vb}{{\bf b}}
\newcommand{\bx}{{\boldsymbol x}}
\newcommand{\by}{{\boldsymbol y}}
\newcommand{\bxero}{{\boldsymbol 0}}
\newcommand{\binf}{\boldsymbol{\infty}}

\title{Low-dimensional optimal stopping with rank-dependent loss}
\author{Alexander Gnedin\and Marcos Costa Santos Carreira}
\date{18 November 2026}

\begin{document}
\maketitle
\begin{abstract}
\noindent
Sequential selection problems with independent observable scores and rank-dependent objectives
are not reducible, in general, to dynamic optimisation over a state space of fixed dimension.
We consider optimal stopping problems with truncated loss structures through the lens of the planar Poisson process,
where self-similarity removes the time variable and the state is reduced to the $m$ lowest scores observed so far, rescaled to the initial time.
For $m=1$ the optimal rule stops at the first record below a hyperbola, which is found by a sliding-window argument.
For $m=2$ we derive the equations for the boundaries of the stopping sets of $1$-records and $2$-records, solve them explicitly on a part of the state space, and express the optimal values by quadratures. 
The value function and the optimal strategy  are  characterised by a static
integro-differential equation. The optimal strategies  are compared with simple low-parameter rules, which are close to optimality.

\end{abstract}

\section{Introduction}

Let $X_1,\ldots,X_n$ be independent random scores sampled from some continuous distribution.
The scores are observed sequentially, with the aim to select one of them
using a nonanticipating stopping rule $\tau$ with values in $[n]:=\{1,\ldots,n\}$.
With each score $X_k$ we associate its {\it final} rank among all $n$ scores,
\begin{equation}\label{f-rank}
R_{k,n}:=\#\{i\in [n]:  X_i\leq X_k\}.
\end{equation}
In `secretary problems' the loss for selecting $X_k$ depends on the scores $(X_1,\ldots,X_n)$
through the ranking permutation $(R_{1,n},\ldots,R_{n,n})$ \cite{Samuels}. This includes in the first turn stopping problems with
 objectives of the form
\begin{equation}\label{q-form}
{\mathbb E}\,q(R_{\tau,n})\to \min,
\end{equation}
where $q(r)$ is a nondecreasing function on the integers; the classical best-choice objective being $q(r)={\bf 1}(r>1)$.  
The nature of a selection problem largely depends on the type of information allowed in stopping rules.

In 
{\it rank-information} problems one considers the family of stopping rules adapted to the sequence $(I_k, ~k\in [n])$ of {\it initial} ranks
\begin{equation}\label{i-rank}
I_{k}:=\#\{i\in [k]:  X_i\leq X_k\},
\end{equation}
where $I_1,\ldots,I_k$ encode the ranking permutation of the first $k$ scores.
Adopting the terminology of extreme-value theory, we say that $X_k$ is an {\it $i$-record} if $I_k=i$, and just a {\it record} if $i=1$.
By exchangeability, the sequence of final ranks is a uniformly random permutation of $[n]$, the initial ranks are independent with $I_k$ uniformly distributed on $[k]$, and $R_{k,n}$ is independent of $I_1,\ldots,I_{k-1}$.
For problems (\ref{q-form}) this implies that the optimal decision at stage $k$ depends on the observed scores through $I_k$, which makes the stopping problems efficiently solvable by backward induction, as the continuation value only depends on $k$ and $n$.
The structural clarity underlies the large-sample theory developed for (\ref{q-form}) 
and related problems in the framework of a continuous-time infinite-arrivals model \cite{GiSamuels},
where the value function is determined from a piecewise differential equation \cite{Mucci}.
Closed-form asymptotic solutions are available, in particular, 
for the problem of choosing one of the $m$ best with $q(r)={\bf 1}(r>m)$  \cite{FS, GM, GZ, Mucci}, and Lindley's 
minimum rank problem with the linear loss $q(r)=r$ \cite{CMRS, GnedinUniform, Mucci}.

In {\it full-information} problems the optimum is sought within the larger family of all stopping rules adapted to $(X_k, ~k\in [n])$. The term `full-information' was coined in \cite{GM} to stress that the observer knows exactly the source distribution (which may be assumed uniform), and can use the score values in the selection strategies.
The informational advantage makes the optimisation problem much more complex, because the global ranking couples the observed and future scores, thus forcing history 
dependence  of the continuation value and decisions.
To date the sole completely understood instance is the best-choice problem, reducible to optimal stopping of the univariate Markovian sequence of records \cite{GFI, BC-SPA, GKS}.
For problems (\ref{q-form}) with general $q$, in particular for the linear loss -- known in this full-information setting as the Robbins problem -- one should not expect a closed-form solution, since the continuation value at stage $k$ depends on all $k$ past scores \cite{Bruss, Swan}, and the full history dependence persists in the limit \cite{GnedinJAP}.
The present paper aims to illustrate how the complexity grows with the dimension, by focusing on truncated loss structures, where only a small number $m$ of the $i$-record processes, $i\le m$, is relevant.
Concretely, we consider two loss structures:
\begin{eqnarray}\label{Adis}
{\rm A}_{n,m}\!:&\!\!\!\!\! \text{the loss for selecting an $i$-record } X_k \text{ is } \min(i,\,m+1)+\#\{k<j\leq n: X_j<X_k\},\\
\label{Bdis}
{\rm B}_{n,m}\!:& q(r)=\mathbf 1(r>m)\ \text{ for the final rank } r.
\end{eqnarray}
In problem ${\rm A}_{n,m}$ the loss is the capped initial rank plus the number of later smaller scores. This objective was employed in \cite{BF} to compute numerical lower bounds for the value in Robbins' problem.
We will show that problem ${\rm A}_{n,1}$ has a simple asymptotic solution, very much in line with the best-choice problem ${\rm B}_{n,1}$,
both are treated here by a novel sliding window approach.
We further consider ${\rm A}_{n,2}$ and ${\rm B}_{n,2}$, which require keeping track of the $1$- and $2$-record processes.
For ${\rm B}_{n,m}$ the best of the explored single-curve rules has the failure probability within $7.5\%$ ($m=2$) and about $20\%$ ($m=3$) of the optimum computed in \cite{GMal}, and 
two- or three-parameter rules reduce the excess to below $1\%$. This suggests that for problems with rank-dependent loss the complexity of the exactly optimal strategies is not matched by the gain over simple low-parameter rules,
 in accord with our findings on  the best-choice problem \cite{Memoryless}.

In \cite{GFI, BC-SPA, GKS} the full-information best-choice problem and some of its versions were embedded in a planar Poisson process (PPP) --
an infinite-arrivals model providing a natural framework for the asymptotic analysis.
In \cite{GnedinJAP} this connection was extended to  problems with objective (\ref{q-form}) where the stopping value remains bounded as $n\to\infty$.
In \cite{Memoryless} we applied the PPP as a tool to compare the optimal best-choice strategy against a memoryless rule, which stops at the first arrival below a precommitted time-dependent threshold.
Here we continue this line by using the self-similarity of the PPP to eliminate the time variable from the integro-differential optimality equation, 
 to reveal the structural patterns of the optimal adaptive threshold, and to evaluate the limits in (\ref{Adis}) and (\ref{Bdis}).

The rest of the paper is organised as follows. In Section~\ref{sec:setting} we connect the discrete-time problems (\ref{Adis}), (\ref{Bdis}) with their PPP counterparts $\mathrm A_m,\mathrm B_m$, show the convergence of the optimal stopping values, and analyse the optimality equation for the value function. The central optimality result, Theorem~\ref{thm:bellman}, 
is proved for more general loss structures;   we need a special argument using
self-similarity and an embedded  Markov chain with a finite number of steps, since the PPP setting, with infinitely many arrivals and a rank-dependent loss, is not covered verbatim by the standard theory.
Proposition~\ref{prop:mon} describes the pathwise behaviour of the adaptive acceptance threshold: we show that this is a nondecreasing process drifting continuously between arrivals, jumping up at some arrival times, and staying flat once it reaches a past score.
Section~\ref{sec:one} presents explicit solutions for $m=1$. Section~\ref{sec:two} describes the geometric boundaries  of the two-dimensional stopping sets, and Sections~\ref{sec:A2} and \ref{sec:B2} give the solutions of $\mathrm A_2$ and $\mathrm B_2$. Section~\ref{sec:thr} explores the performance of simple strategies, especially for $\mathrm B_m$. Section~\ref{sec:num} presents numerical evaluation of optimal and suboptimal stopping strategies.

\section{The planar Poisson process framework}\label{sec:setting}

\subsection{Scores, ranks and losses}
Let $\Pi$ be a unit rate PPP on $[0,1]\times[0,\infty)$. We view $\Pi$ as a random scatter of points, where an atom $(t,x)\in\Pi$ is interpreted as a score $x$ observed at time $t$, smaller scores being preferred. The natural filtration is generated by the restrictions of $\Pi$ to $[0,t]\times[0,\infty)$, $t\in[0,1]$.
An admissible  selection strategy (sometimes referred to as rule)  is a stopping time $\tau$ with respect to this filtration, which is either the arrival time of an atom or $\tau=1$ (no selection); the latter happens with positive probability for the optimal strategies. 
See \cite{Rusch} for measure-theoretic aspects of optimal stopping for planar point processes,
with objectives depending on $(t,x)$ alone.

For an atom $(t,x)$ we define its {\it final} rank $r$ to be one plus the number of atoms of $\Pi$ with scores smaller than $x$, and its {\it initial} rank $i$ to be one plus the number of earlier atoms with scores smaller than $x$. When the atom arrives,  $x$ and $i$ are observed, and the conditional distribution of $r-i$ is Poisson$(\xi)$, with the mean equal to the area  of the south-east window   $(t,1]\times(0,x)$, where
$$\xi:=\bar{t}x,\qquad \bar t:=1-t .$$
We call $\xi$ the {\it box coordinate} of the atom.

For $m=1,2,\dots$ we consider two types of the {\it adapted} loss from stopping at an atom with initial rank $i$ and box coordinate $\xi$, given the observed information:
\begin{align}
&\mathrm A_m:&& Q(\xi,i)=\min(i,m+1)+\xi, \label{QA}\\
&\mathrm B_m:&& Q(\xi,i)=1-\sum_{j=0}^{m-i} e^{-\xi}\frac{\xi^{j}}{j!}, \label{QB}
\end{align}
where the sum in (\ref{QB}) is empty for $i>m$. 
The PPP problems with these loss structures are infinite-arrivals counterparts of the full-information problems 
(\ref{Adis}), (\ref{Bdis})
from the Introduction.
In both cases $Q(\xi,i)$ is nondecreasing in $\xi$ and in $i$. 
For the event $\{\tau=1\}$ we assume the loss
$$
Q_\infty=m+1~~~{\rm in~~}{\rm A}_m,\qquad Q_\infty=1~~~{\rm in~~}{\rm B}_m .
$$
This convention is the limit form of the terminal loss in the finite-$n$ problems (\ref{Adis}),(\ref{Bdis}), where the observer must stop at stage $n$: by exchangeability, the last score has rank of order $n$ with high probability, and the loss for being forced to accept it tends to $Q_\infty$. Clearly, it can only be optimal to stop at observations with initial rank $i\leq m$, since otherwise the loss is at least $Q_\infty$; this may lead to missing all opportunities and ending at $t=1$.
The minimal achievable expected losses are denoted ${\bf a}_m$ and ${\bf b}_m$ for problems ${\rm A}_m$ and ${\rm B}_m$, respectively.
For ${\rm B}_m$ the loss is the conditional probability of {\it failure}, which is the choice of an atom of final rank above $m$; the complementary probability $1-\vb_m$ is the probability of successful choice of one of the $m$ best scores.

Here is the link to the discrete-time fixed-$n$ problems. Split the time range $[0,1]$ in $n$ equal slots $((k-1)/n,\,k/n],~k\in[n]$, and suppose that the arrival times of the $n$ slot-wise minima
$M_1,\dots,M_n$
 belong to the prior information of the observer. Then the selection is reduced to an $n$-step problem for $M_1,\dots,M_n$, which are independent and share the same exponential distribution. 
With the further advantage that only unobserved minima contribute to the final rank of the selected observation, the problem becomes equivalent to ${\rm A}_{n,m}$, respectively ${\rm B}_{n,m}$. This shows that the Poisson values ${\bf a}_m,{\bf b}_m$ are upper bounds for the values of the discrete problems. 
In problem (\ref{Bdis}) the values increase with $n$ to ${\bf b}_m$ \cite{GnedinJAP}. For problem (\ref{Adis}) the same holds, but requires a separate  proof, since the loss is not a function of the final rank.

\begin{proposition}\label{prop:conv}
Let $a_{n,m}$ be the value of problem {\rm(\ref{Adis})}. For every $m\ge1$ the sequence $a_{n,m}$, $n\ge m$, is nondecreasing and converges to ${\bf a}_m$.
\end{proposition}
\begin{proof}
The loss for selecting a score is $\min(i,m+1)+N$, where $i$ is its initial rank and $N$ is the number of later scores below it. It depends on the scores only through their order, hence
 the common continuous distribution of the scores does not matter for the stopping problem.

 Monotonicity in $n$ follows by a familiar coupling argument. Take $n+1$ independent uniform scores, $n\ge m$, and let the observer be told in advance the value $u$ of the largest one; the advantage can only lower the minimal expected loss. If $u$ is the $k$-th score, selecting it costs $\min(k,m+1)+(n+1-k)\ge m+1$, while the last of the other scores
costs at most $m+1$,
so $u$ may be declined. 
Discarding $u$ does not impact ranking of the other scores,
and given $u$ the other scores are independent and uniform on $[0,u]$. Thus the problem is reduced to $n$ scores, and $a_{n,m}\le a_{n+1,m}$.

 Upper bound follows by the slot argument above: $a_{n,m}\le{\bf a}_m$.  The rest of the proof deals with a lower bound.

Fix $n$, a level $A>0$ and $\delta\in(0,1)$, and let $M_1,\dots,M_n$ be the minima of $\Pi$ over the slots $((k-1)/n,\,k/n]$, as above; they are independent and $\P[M_k>x]=e^{-x/n}$.
Let $\rho$ be the optimal strategy for the $n$-step problem with the scores $M_1,\dots,M_n$,
i.e. $\rho$ is adapted to  $\mathcal G_k:=\sigma(M_1,\dots,M_k), k\in [n]$, and incurs some loss $L_\rho$ satisfying $\E L_\rho=a_{n,m}$.
At stage $k$ the strategy $\rho$ accepts any score below
some  threshold depending on history $M_1,\dots,M_{k-1}$.  The latter follows 
from the opposite directions of monotonicity of the expected loss of stopping at score $x$, given by 
$\min(i,m+1)+(n-k)(1-e^{-x/n})$ with $i$ the initial rank of $x$, 
 and  the value function, which  is nonincreasing in the state (the coupling of Lemma~\ref{lem:mon}(i)). 
Consider a stopping time $\tau$ for $\Pi$ defined as follows: $\tau$ stops at an atom $(t,x)$ of the $k$-th slot if and only if $x<A$, no earlier atom of this slot has score below $A$, and $\rho$ applied to the sequence $M_1,\dots,M_{k-1},x$ stops at stage $k$ with score $x$. At any time $t\in[0,1]$ the minima of the earlier slots are known, so $\tau$ is a legitimate stopping time for $\Pi$. Let $L_\tau$ be its loss, equal to $m+1$ if no atom is selected, and let $N_A\stackrel{d}{=}{\rm Poisson}(A)$ count  the number of atoms of $\Pi$ with score below $A$.

Let $G$ be the event that every slot contains at most one atom with score below $A$, and let $D=\{M_\rho\ge A\}$.
Suppose $\tau$ stops at an atom $x$ of the $k$-th slot. Then $x\ge M_k$, and $\rho$ applied to $M_1,\dots,M_{k-1},x$ stops at the stage $k$ and not before;  
then  $M_k\le x$ belongs to the acceptance interval for the same history, hence $\rho=k$. Thus $\tau$ stops only in the slot $\rho$, at the first atom of this slot with score below $A$ (which is not later than the atom $M_\rho$ in time), or makes no selection.
On $G\setminus D$ we have $L_\tau=L_\rho$, since
the slot $\rho$ contains exactly one atom below $A$, which has the score $M_\rho$ selected
by $\tau$; every atom below $M_\rho$ is the minimum of its slot, so the initial rank of 
$M_\rho$ and the number of later atoms below it are the same in $\Pi$ as among the slot minima.
On $G\cap D$ the slot $\rho$ has no atom below $A$, thus $\tau$ makes no selection and $L_\tau=m+1$.
On the complement of $G$ the selected atom may differ from the slot minimum, but still $L_\tau\le m+1+N_A$. Since $L_\rho\ge0$, we obtain
\begin{equation}\label{eq:Ltau}
\begin{aligned}\E L_\tau&\le\E[ L_\rho\mathbf 1(D^c)]+(m+1)\P[D]+\E[(m+1+N_A)\mathbf 1(G^c)]\\
&=a_{n,m}+\E[(m+1-L_\rho)\mathbf 1(D)]+\E[(m+1+N_A)\mathbf 1(G^c)].\end{aligned}
\end{equation}

{\it The term with $D$.} Since $D\in\mathcal G_\rho$ and, conditionally on $\mathcal G_\rho$, the later minima are independent of the past, on $D$ we have
$$\E[L_\rho\,|\,\mathcal G_\rho]=\min(i_\rho,m+1)+(n-\rho)\big(1-e^{-M_\rho/n}\big)\ge\min(i_\rho,m+1)+(n-\rho)\big(1-e^{-A/n}\big),$$
where $i_\rho$ is the initial rank of $M_\rho$. If $\rho\le(1-\delta)n$, the right side is at least $\delta n(1-e^{-A/n})$, which tends to $\delta A$ as $n\to\infty$, so it exceeds $m+1$ for large $n$ when $A\ge2(m+1)/\delta$. If $\rho>(1-\delta)n$, then on $D$ every minimum $M_k$ with $k\le(1-\delta)n$ and $M_k<A$ is below $M_\rho$ and arrives earlier, thus $i_\rho\ge m+1$ unless fewer than $m$ such minima exist. Let $H$ be the event that fewer than $m$ of the minima $M_k$, $k\le(1-\delta)n$, are below $A$. Then, for large $n$, $(m+1-\E[L_\rho|\mathcal G_\rho])\mathbf 1(D)\le(m+1)\mathbf 1(H)$, which yields
$$\E[(m+1-L_\rho)\mathbf 1(D)]\le(m+1)\P[H].$$
The number of  minima $M_k$, $k\le(1-\delta)n$, which fall below $A$ is binomial with success probability $1-e^{-A/n}$, hence $\P[H]\to\P[N\le m-1]$ as $n\to\infty$, where $N$ is Poisson with mean $(1-\delta)A$.

{\it The term with $G^c$.} Let $N_k$ be the number of atoms with score below $A$ in the $k$-th slot; these are independent Poisson  with mean $A/n$, and $N_A=\sum_kN_k$, $G^c=\bigcup_k\{N_k\ge2\}$. Therefore $\P[G^c]\le n\,(A/n)^2/2=A^2/(2n)$, and
\begin{align*}\E[N_A\mathbf 1(G^c)]&\le\sum_{j\ne k}\E N_j\,\P[N_k\ge2]+\sum_k\E [N_k(N_k-1)]\\&\le n^2\,\frac An\,\frac{A^2}{2n^2}+n\,\frac{A^2}{n^2}=\frac{A^3}{2n}+\frac{A^2}{n},\end{align*}
since $N\le N(N-1)$ for $N\ge2$ for the diagonal terms $j=k$. Therefore for the term in focus
$\E[(m+1+N_A)\mathbf 1(G^c)]\to0$ as $n\to\infty$.

Finally, given the past, the number of later atoms below a selected atom is Poisson with mean equal to its box coordinate, so $\E L_\tau$ is the expected loss (\ref{QA}) of the strategy $\tau$, and ${\bf a}_m\le\E L_\tau$. Passing to the limit in (\ref{eq:Ltau}) we obtain for every $A\ge2(m+1)/\delta$
$${\bf a}_m\ \le\ \lim_{n\to\infty}a_{n,m}+(m+1)\,\P[N\le m-1],$$
and letting $A\to\infty$ gives ${\bf a}_m\le\lim_na_{n,m}$ (the limit exists by monotonicity). Taken together with the upper bound, this implies the desired convergence $a_{n,m}\uparrow{\bf a}_m$.
\end{proof}

\subsection{Self-similarity and the state}
For every $t>0$ the number of atoms arrived before time $t$ is infinite. However, since the ranks above $m$ can be ignored, the past data relevant for the future decisions is the ordered vector $\by=(y_1,\dots,y_m)$, where $y_i$ is the $i$-th lowest score observed before $t$. We call $\by$ the {\it history}; this statistic takes the role of the state 
in a Markovian stopping problem. 
Along the path of PPP, the history is a monotonic jump Markov process; the marginals $y_i$ are the (minimal) $i$-record processes,
coming down from infinity at $t=0+$  \cite{GnedinCorners}. 
In the sequel under time-$0$ history we shall mean a finite state of this $m$-dimensional process; the  stopping problem 
initialised this way
can be interpreted in terms of Palm conditioning of the PPP on $m$ arrivals at time $0$.

\begin{lemma}\label{lem:self}
The decision problem at time $t$ with the history $\by$ is equivalent to the problem at time $0$ with the initial history $\bar t\by=(\bar ty_1,\ldots,\bar ty_m)$.
\end{lemma}
\begin{proof}
Given the observations up to time $t$, the restriction of $\Pi$ to $(t,1]\times[0,\infty)$ is a PPP. The map $(s,x)\mapsto((s-t)/\bar t,\ \bar tx)$ preserves the Lebesgue measure and sends this restriction to a PPP on $[0,1]\times[0,\infty)$. It sends the history $\by$ to $\bar t\by$, preserves all ranks, and the box coordinate $\bar tx$ of an atom $(t,x)$ is the same in both settings.
\end{proof}

Thus it suffices to consider the stopping--continuation dilemma at time $0$, for the generic state $\bx=(x_1,\dots,x_m)$, where $x_i$ is the $i$-th lowest of the scores in { box coordinates}, $x_i=\bar{t}y_i$. To avoid confusion, we keep the notation $\by$ for the history in the original PPP coordinates at a given time $t>0$.
Let $\mathcal S=\{\bx:\ 0\le x_1\le\dots\le x_m<\infty\}$ be the set of finite states. For $\bx\in\mathcal S$ and $\xi<x_m$ we write $\bx\cup\xi$ for the vector of the $m$ lowest values among $x_1,\dots,x_m,\xi$; this is the updated state upon arrival at the box coordinate $\xi$. The initial rank of such an arrival is $i(\bx,\xi):=1+\#\{j\in[m]: x_j<\xi\}$. An atom with box coordinate above $x_m$ has initial rank larger than $m$; it is never worth accepting (its loss is at least $Q_\infty$) and does not change the history.

\subsection{The value function and the optimal strategy}\label{sec:verif}
The optimality results below apply to a general loss $Q(\xi,i)\ge0$,   provided it is nondecreasing in $\xi\geq0$ and  $i\geq 1$ and satisfies $Q(\xi,i)\ge Q_\infty$ for $i>m$, where $Q_\infty<\infty$ is the loss of no selection, and also $Q(0+,1)<Q_\infty$. These properties hold for the losses (\ref{QA}) and (\ref{QB}). We  require further that
\begin{equation}\label{eq:H}
\liminf_{\xi\to\infty}Q(\xi,1)\ \ge\ Q_\infty ,
\end{equation}
which holds in both cases: $Q(\xi,1)=1+\xi$ in $\mathrm A_m$, and $Q(\xi,1)=1-e^{-\xi}\sum_{j<m}\xi^j/j!\to1$ in $\mathrm B_m$. Condition (\ref{eq:H}) is used only for the attainment of the value for the empty history, see Remark~\ref{rem:H}. Lemma~\ref{lem:mon} below uses further properties of $\mathrm A_m$ and $\mathrm B_m$.

The expected loss of an admissible strategy $\tau$ in the problem which starts at time $0$ with the history $\bx$ is
$$L_{\bx}(\tau)=\E\big[Q(\xi_\tau,i_\tau)\mathbf 1\{\tau<1\}+Q_\infty\mathbf 1\{\tau=1\}\big],$$
where $\xi_\tau,i_\tau$ are the box coordinate and the initial rank (counted together with the history $\boldsymbol x$) of the selected atom; for $\mathrm A_m,\mathrm B_m$ this is the expected loss in the original problem, by the conditioning in Section~\ref{sec:setting}. The history is a state $\bx$ with components in $[0,\infty]$: a finite state $\bx\in\mathcal S$, the empty history $\binf$, or a state with some components equal to $\infty$ (fewer than $m$ scores observed). The {\it value function} is $v(\bx)=\inf_\tau L_{\bx}(\tau)$; thus $v(\binf)$ is ${\bf a}_m$, respectively ${\bf b}_m$ in the corresponding stopping problems. For bounded Borel 
functions $w$ on $\mathcal S$ define the operator
$$P[w](\bx):=\int_0^{x_m}\min\Big\{Q\big(\xi,i(\bx,\xi)\big),\ w(\bx\cup\xi)\Big\}\,\dd\xi .$$

The map $(t,x)\mapsto(u,\xi)=\big(-\log(1-t),\,(1-t)x\big)$ preserves the Lebesgue measure, and transforms  $\Pi$ into a unit rate PPP $\tilde\Pi$ on $[0,\infty)\times[0,\infty)$; the time $t=1$ corresponds to $u=\infty$. In the new time variable $u$ the state in box coordinates, $Z_u=\bar t\,\by_t$, is a piecewise deterministic Markov process: 
between arrivals it drifts from some state $Z_u$ along  
the ray $Z_{u'}=e^{-(u'-u)}Z_u$ toward the origin.
Arrivals with box coordinates $\xi<Z_{u,m}$ that change the state to $Z\cup\xi$ occur at rate $Z_{u,m}$ with $\xi$ uniform on 
$[0,Z_{u,m}]$. (The map is the self-similarity of Lemma~\ref{lem:self} written for all times at once.) The loss $Q(\xi,i)$ depends on an atom only through $\xi$ and the initial rank, and the stopping problem in these coordinates does not involve time.

A constant threshold $b$ for the box coordinate is the hyperbola $x=b/\bar t$ in the original coordinates.
For $m=1$ the state is the box coordinate of  the running minimum,
which makes the stopping problem one-dimensional and implies that 
in the original coordinates the optimal threshold for record scores is 
a hyperbola (Section~\ref{sec:one}).

\begin{remark}\label{rem:hyperbola}
Hyperbolas appear also for the simple strategies of Section~\ref{sec:thr}, which do not use the history.
For $m\ge2$ the threshold depends on the history through the state $\bar t\by$, which shrinks between arrivals while the past scores stay fixed, so the threshold is not a hyperbola (Section~\ref{sec:two}).
The same scaling gives a hyperbola for the problem of  minimising the expected score (Moser's problem in the PPP \cite{GnedinUniform}), where the optimal threshold is $2/\bar t$.
In discrete time there is no self-similarity, and hyperbolic thresholds appear only asymptotically \cite{GM, GFI, BC-SPA}.
\end{remark}

Call an atom $(t,x)$ {\it low} if its initial rank is at most $m$, that is, if $x<y_m$ for the history $\by=\by_{t-}$ just before $t$, or equivalently if its box coordinate $\xi=\bar tx$ is below the $m$-th component of the state $Z=\bar t\,\by$ just before $t$. A low atom changes the state by a jump when it arrives, and it may be accepted or declined. An atom which is not low does not change the state and is declined.

Consider the problem which starts at time $0$ with the history $\bx\in\mathcal S$; at time $0$ the original and the box coordinates coincide. The history $\by_t$ in the original coordinates can only decrease, so $y_m\le x_m$ at all times. Hence every low atom has score below $x_m$, and the number $N$ of low atoms is stochastically dominated by the Poisson variable with mean $x_m$, 
\begin{equation}\label{eq:Nrel}
\P[N>k]\le\P[\mathrm{Poisson}(x_m)>k]\to0\quad(k\to\infty),\ \ \text{uniformly for }x_m\le a.
\end{equation}

The selection problem is thus an optimal stopping problem for a Markov chain with an a.s.\ finite number of steps: the chain is the sequence of states $(i,Z)$ right after the successive low atoms,
with loss $Q(\xi,i)$ being a function of the state, where $\xi=Z_i$.
This makes applicable the standard discrete-time Markovian stopping theory (see \cite{RSC}, Section 5.1).

By the first-arrival decomposition to follow, the first low atom after the time $0$, in the problem with the state $\bx$, arrives at the time $r$ in the variable $r=e^{-u}\in(0,1)$ with the density $x_me^{-x_m(1-r)}\dd r$, its box coordinate is uniform on $[0,rx_m]$ and independent of $r$, and the state after the atom is $r\bx\cup\xi$; with probability $e^{-x_m}$ there is no low atom at all. For a Borel set $A\subset\{(\bx,\xi):\xi<x_m\}$ (referred to as {\it acceptance set} for the arriving score) and a bounded Borel $w$ define
\begin{equation}\label{eq:rayA}
(T_Aw)(\bx)=e^{-x_m}\Big(Q_\infty+\int_0^1e^{x_mr}\,\frac{P_A[w](r\bx)}r\,\dd r\Big),
\end{equation}
$$P_A[w](\bx)=\int_0^{x_m}\Big(\mathbf 1_A\,Q+\mathbf 1_{A^c}\,w(\bx\cup\xi)\Big)\dd\xi,\qquad Q=Q\big(\xi,i(\bx,\xi)\big),$$
and $T=T_{A_w}$ for $A_w=\{(\bx,\xi):Q(\xi,i(\bx,\xi))\le w(\bx\cup\xi)\}$, so that $P_{A_w}[w]=P[w]$ and
\begin{equation}\label{eq:ray}
(Tw)(\bx)=e^{-x_m}\Big(Q_\infty+\int_0^1e^{x_mr}\,\frac{P[w](r\bx)}r\,\dd r\Big).
\end{equation}
If $0\le w\le Q_\infty$ and $Q\le Q_\infty$ on $A$ (in particular for $A=A_w$), then $0\le P_A[w](r\bx)\le rx_mQ_\infty$, so the integral converges and $e^{-x_m}Q_\infty\le T_Aw\le Q_\infty$.

The next result makes the dynamic programming equation of this chain explicit, and adds the uniqueness of its solution, the optimal strategy, and the case of the empty history.

\begin{theorem}\label{thm:bellman}
{\rm(i)} The value function $v$ is continuous on $\mathcal S$ with $0\le v\le Q_\infty$ and $v(\bxero)=Q_\infty$, and along every ray the function $s\mapsto v(s\bx)$ is continuously differentiable on $(0,\infty)$, with
\begin{equation}\label{eq:bellman}
\sum_{i=1}^m x_i\frac{\partial v(\bx)}{\partial x_i}+x_m\,v(\bx)=P[v](\bx),
\end{equation}
where $\sum_ix_i\partial_iv(\bx)$ is the derivative of $s\mapsto v(s\bx)$ at $s=1$. Equivalently, $v=Tv$.

{\rm(ii)} The solution is unique: every bounded solution of $w=Tw$ is equal to $v$. In particular, a function $w$ which is bounded, continuous and continuously differentiable along the rays, satisfies {\rm(\ref{eq:bellman})} and $w(s\bx)\to Q_\infty$ as $s\to0$, is the value function.

{\rm(iii)} The strategy $\tau^*$ which, in the problem with the initial state $\bx\in\mathcal S$, selects the first low atom with the box coordinate $\xi$ satisfying
\begin{equation}\label{eq:balance}
Q\big(\xi,i(Z,\xi)\big)\le v(Z\cup\xi),
\end{equation}
where $Z$ is the state immediately before the arrival of the atom, is optimal: $L_{\bx}(\tau^*)=v(\bx)$.

{\rm(iv)} The value for the empty history is $v(\binf)=\lim_{c\to\infty}v(c,\dots,c)=\inf_{\bx} v(\bx)$. The strategy $\tau^*$, defined by {\rm(\ref{eq:balance})} with $Z$ the current history (the missing components being $\infty$), is optimal for the empty history whenever the stopping time is not degenerate at $0$, that is, $\P[\tau^*>0]=1$; this holds under {\rm(\ref{eq:H})}.
\end{theorem}
\begin{proof}
{\it Step 1} (the contraction estimate). Let $w,w'$ be bounded Borel functions and 
$$\delta(q)=\sup\{|w(\by)-w'(\by)|:y_m\le q\}.$$ 
For every acceptance set $A$, $|P_A[w]-P_A[w']|(r\bx)\le rx_m\,\delta(rx_m)$, and the same holds for $P$ since $|\min\{Q,w\}-\min\{Q,w'\}|\le|w-w'|$. Hence, by (\ref{eq:ray}) and $e^{-x_m(1-r)}\le1$,
\begin{equation}\label{eq:gron}
|T_Aw(\bx)-T_Aw'(\bx)|\vee |Tw(\bx)-Tw'(\bx)|\ \le\ \int_0^1x_m\,\delta(rx_m)\,\dd r=\int_0^{x_m}\delta(q)\,\dd q.
\end{equation}
If $w=T_Aw$ and $w'=T_Aw'$ (or $w=Tw$, $w'=Tw'$), then $\delta(a)\le\int_0^a\delta(q)\dd q$ for all $a$, with $\delta$ bounded and nondecreasing, so by iteration $\delta(a)\le\sup\delta\cdot a^n/n!\to0$ and $w=w'$. This is the uniqueness of the bounded Borel fixed point of $T_A$ and of $T$.

{\it Step 2} (iteration). Let $w_0\equiv Q_\infty$, $w_{k+1}=Tw_k$. The operator $T$ is monotone, and $w_1\le w_0$, hence $w_{k+1}\le w_k$, and $w_k\downarrow\hat v$ with $e^{-x_m}Q_\infty\le\hat v\le Q_\infty$. Each $w_k$ is continuous: for bounded continuous $w$ the map $\bx\mapsto P[w](\bx)$ is continuous, because for every $\xi\notin\{x_1,\dots,x_m\}$ the integrand and the upper limit of integration depend continuously on $\bx$, and the integrand is bounded by $Q_\infty$ (by dominated convergence); then $Tw$ is also continuous by the same token. With $\delta_k(a)=\sup_{x_m\le a}|w_{k+1}-w_k|$, (\ref{eq:gron}) gives $\delta_{k+1}(a)\le\int_0^a\delta_k(x){\rm d}x$, hence $\delta_k(a)\le Q_\infty a^k/k!$, so $w_k\to\hat v$ uniformly on $\{x_m\le a\}$ for every $a$, $\hat v$ is continuous, satisfies $\hat v=T\hat v$ by monotone convergence, and $\hat v(\bxero)=Q_\infty$.

{\it Step 3} (the value of the problem with $k$ decision epochs). 
The decision epochs are the arrivals of low atoms; and the sequence of renewed histories is a Markov chain terminating after finitely many steps (no further low atoms), with the terminal loss
$Q_\infty$.
Let $\mathcal R_k$ be the class of admissible strategies which may select only from the first $k$ low atoms.  Optimisation over $\mathcal R_k$
 is a finite-horizon problem, solved by backward induction \cite{RSC}: conditionally on the first low atom,  by Lemma~\ref{lem:self} and the Poisson property the future is the problem started at the new state $r\bx\cup\xi$, so the value of the $k$-epoch problem is $w_k=Tw_{k-1}$, $w_0=Q_\infty$, which are the functions of Step 2. Now let $\tau$ be an arbitrary strategy and let $\tau^{(k)}\in\mathcal R_k$ make the same selection if it is a low atom among the first $k$, and no selection otherwise. Selecting an atom of initial rank above $m$ costs at least $Q_\infty$, so the loss changes only on the event that $\tau$ selects a later atom, and $L_{\bx}(\tau^{(k)})\le L_{\bx}(\tau)+Q_\infty\P[N>k]$. By (\ref{eq:Nrel}) and $w_k\ge\hat v$,
$$L_{\bx}(\tau)\ge w_k(\bx)-Q_\infty\P[N>k]\ge\hat v(\bx)-Q_\infty\P[N>k]\to\hat v(\bx),$$
whence $v\ge\hat v$.

{\it Step 4} (the strategy $\tau^*$). Let $A^*=A_{\hat v}$ be the acceptance set of (\ref{eq:balance}) with $\hat v$ in place of $v$, and $W(\bx)=L_{\bx}(\tau^*)\in[0,Q_\infty]$ (a selected atom has $Q\le\hat v\le Q_\infty$). The first-arrival decomposition gives $W=T_{A^*}W$, and $\hat v=T\hat v=T_{A^*}\hat v$, since the minimum is attained by $Q$ on $A^*$. By the uniqueness in Step 1, $W=\hat v$. With Step 3, $\hat v=W\ge v\ge\hat v$. Hence $v=\hat v$, and $\tau^*$ is optimal. This proves (iii), the uniqueness in (ii) for $T$, and $v=Tv$ in (i).

{\it Step 5} (the differential form). By Step 2 and $v=\hat v$, for fixed $\bx$ the function $f(s)=v(s\bx)$ satisfies $$f(s)=e^{-sx_m}\big(Q_\infty+\int_0^se^{x_m\rho}P[v](\rho\bx)\rho^{-1}\dd\rho\big)$$ with a continuous integrand bounded by $x_mQ_\infty e^{x_m\rho}$, so $f\in C^1(0,\infty)$ and $sf'(s)+sx_mf(s)=P[v](s\bx)$, which is (\ref{eq:bellman}) at $s\bx$. Conversely, let $w$ be as in the second statement of (ii). Then $(e^{sx_m}w(s\bx))'=e^{sx_m}P[w](s\bx)/s$, and integrating over $[\varepsilon,1]$ and letting $\varepsilon\to0$, with $w(\varepsilon\bx)\to Q_\infty$, gives $w=Tw$, so $w=v$ by Step 1. The same computation shows that, for a bounded $w$ which is continuously differentiable along the rays, $w=Tw$ is equivalent to (\ref{eq:bellman}) together with $w(s\bx)\to Q_\infty$ as $s\to0$.

{\it Step 6} (the empty history). Let $t>0$ and let $Z_t:=\bar t\by_t$ be the state at time $t$ (indexed in this step by the original time) in the problem with the empty initial history; $Z_t$ is a.s.\ finite and  tends to $\binf$ as $t\to0$.
Conditionally on the information at  time $t$ the stopping problem is the one with the initial history  $Z_t$ (Lemma~\ref{lem:self}). 
 By Lemma~\ref{lem:mon}(i), which is proved by coupling and does not depend on the theorem, $\E v(Z_t)\downarrow v_\infty:=\lim_cv(c,\dots,c)=\inf v$ as $t\downarrow0$.

To argue
{\it $v(\binf)=v_\infty$} we note first that the strategy which does nothing before $t$ and then follows $\tau^*$ has the loss $\E v(Z_t)$ by (iii), so $v(\binf)\le v_\infty$. Conversely, for every strategy $\tau$, conditioning on the information at the time $t$,
$$L_{\binf}(\tau)\ \ge\ \E\big[v(Z_t)\mathbf 1\{\tau>t\}\big]\ \ge\ \E v(Z_t)-Q_\infty\P[\tau\le t]\ \ge\ v_\infty-Q_\infty\P[\tau\le t],$$
and $\P[\tau\le t]\to0$ as $t\to0$, because $\tau>0$ a.s. So $L_{\binf}(\tau)\ge v_\infty$. This does not use (\ref{eq:H}).

Next, we show attainment of the optimum. The strategy $\tau^*$ for the empty history is the first atom satisfying (\ref{eq:balance}), where the components of the state not yet filled are $\infty$; after any time $t>0$ there are only finitely many low atoms, so acceptable atoms can accumulate only at the time $0$. Suppose $\P[\tau^*>0]=1$, that is, $\tau^*$ is not degenerate at $0$. Then $\tau^*$ is a well defined stopping time, and by (iii) applied at the time $t$,
$$L_{\binf}(\tau^*)=\E\big[Q\mathbf 1\{\tau^*\le t\}\big]+\E\big[v(Z_t)\mathbf 1\{\tau^*>t\}\big].$$
Since $\P[\tau^*\le t]\to0$ the dominated convergence gives
 $L_{\binf}(\tau^*)=v_\infty=v(\binf)$.

It remains to see that (\ref{eq:H}) implies $\P[\tau^*>0]=1$. First, $v(\binf)<Q_\infty$: the strategy which stops at the first atom with score below $b$ stops with probability $1-e^{-b}$ at an atom of initial rank $1$ with $\xi\le b$, so its loss is at most $Q_\infty-(1-e^{-b})(Q_\infty-Q(b,1))<Q_\infty$ for small $b$, because $Q(0+,1)<Q_\infty$. By (\ref{eq:H}) and $v(c,\dots,c)\to v(\binf)<Q_\infty$ there is $c_0$ with $Q(\xi,1)>v(\xi,\dots,\xi)$ for all $\xi\ge c_0$. Let $s_0\in(0,\frac12)$ and let $E$ be the event that there is no atom with score below $2c_0$ before $s_0$; $\P[E^c]\le2c_0s_0$. On $E$ every atom before $s_0$ has box coordinate at least $c_0$, and the state after its arrival is at least $(c_0,\dots,c_0)$ componentwise, so by Lemma~\ref{lem:mon}(i) $v(\text{state})\le v(c_0,\dots,c_0)<Q(c_0,1)\le Q(\xi,i)$: no atom before $s_0$ is acceptable. Hence $\P[\tau^*\le s_0]\le2c_0s_0\to0$.
\end{proof}

\begin{remark}\label{rem:H}
The boundary condition is essential: the left side of (\ref{eq:bellman}) vanishes at $\bxero$ regardless of the value there, and the form $v=Tv$ has the boundary condition built in, as the term $Q_\infty e^{-x_m}$ of no low atom. No smoothness in $\bx$ is needed, only along the rays; in fact $v$ is $C^1$ in the interior for $m=2$ (Proposition~\ref{prop:reg}). The strategy $\tau^*$ may terminate without selection.
Condition (\ref{eq:H}) is needed only for the attainment of the value $v(\binf)$ for the empty history, to exclude acceptance of the rank-$1$ atoms with huge scores, which accumulate at $t=0$; the value itself is $\lim_cv(c,\dots,c)$ without it. Without (\ref{eq:H}) the value need not be attained: for $m=1$, $Q(\xi,i)=\mathbf 1\{i\ge2\}$ and $Q_\infty=1$ we have $v(\binf)=0$ (wait until time $s_0$ and accept the next record, which is missed with probability $s_0$) but no strategy achieves this; here every record is acceptable and $\tau^*$ is degenerate at $0$.
\end{remark}

Monotonicity of the value function is shown next.
\begin{lemma}\label{lem:mon}
\begin{itemize}
\item[{\rm(i)}] For every adapted  loss structure $Q$ as above, $v$ is nonincreasing in every component: if $\bx'\le\bx$ componentwise,  then $v(\bx)\le v(\bx')$, where 
some of the components of $\bx$ may be  $\infty$.
\item[{\rm(ii)}] For the losses {\rm(\ref{QA})}, {\rm(\ref{QB})}, $v$ is strictly decreasing in every component $x_1,\ldots,x_m$, and $$v(\binf)< v(\bx)< v(\bxero)=Q_\infty\qquad(\bx\in\mathcal S,\ \bx\ne\bxero).$$
\item[{\rm(iii)}] $v(\bx)\to v(\binf)$ as $x_1\to\infty$.
\end{itemize}
\end{lemma}
\begin{proof}
(i) Let $\tau'$ be a strategy in the problem with the history $\bx'$. It is a function of the future atoms and of $\bx'$, so it can be applied with the history $\bx$. For every realisation the initial rank of the selected atom is not larger under $\bx$ than under $\bx'$, and $Q$ is nondecreasing in the rank, so $L_{\bx}(\tau')\le L_{\bx'}(\tau')$; taking the infimum over $\tau'$ gives $v(\bx)\le v(\bx')$. This proof does not use Theorem~\ref{thm:bellman}, which uses (i) in Step 6.

(ii) For $\bx\ne\bxero$ we have $x_m>0$ and $v(\bx)<Q_\infty$: the strategy which accepts the first low atom with box coordinate below $\varepsilon$ stops with positive probability, at an atom of initial rank at most $m$, and its loss is then at most $Q(\varepsilon,m)<Q_\infty$ for small $\varepsilon$.

For the strict inequalities let $\bx'\le\bx$ componentwise and $\bx'\ne\bx$, where $\bx'\in\mathcal S$ and either $\bx\in\mathcal S$ or $\bx=\binf$; we show that $v(\bx)<v(\bx')$. Let $\tau'$ be the strategy $\tau^*$ of Theorem~\ref{thm:bellman}(iii) for the history $\bx'$; it selects only low atoms, that is, atoms of initial rank at most $m$ under $\bx'$. Since $\tau'$ is a function of the future atoms and of $\bx'$, it can be applied with the history $\bx$, and by the coupling of (i), $L_{\bx}(\tau')\le L_{\bx'}(\tau')$. Suppose that $v(\bx)=v(\bx')$. Then
$$L_{\bx}(\tau')\ \ge\ v(\bx)\ =\ v(\bx')\ =\ L_{\bx'}(\tau')\ \ge\ L_{\bx}(\tau'),$$
so $\tau'$ is optimal for $\bx$ as well, and the two losses coincide. (In general, it is not true that a strategy loses strictly when a component is lowered: a strategy which stops at the first atom with score below $\varepsilon<x'_1$ has the same loss under both histories. So the optimality of $\tau'$ has to be used.)

Let $j$ be an index with $x'_j<x_j$ ($j=1$ if $\bx=\binf$), let $c:=x_m$ if $\bx\in\mathcal S$ and let $c>x'_m$ be a finite level if $\bx=\binf$, and $a:=\min(x_j,c)>x'_j$. Fix $\delta\in(0,1)$, to be chosen below, and let $F$ be the event that $\Pi$ has no atom with score below $c$ in $[0,1-\delta]$, and that the first atom with score below $c$ in $(1-\delta,1]$ has a score $z\in(x'_j,a)$. Then $\P[F]>0$, and on $F$ the strategy $\tau'$ has not stopped before this atom, since an atom selected by $\tau'$ has score below $x'_m\le c$. Let $i<i'$ be the initial ranks of the atom under $\bx$ and under $\bx'$ (the rank is larger under $\bx'$ because $x'_j<z<x_j$), so that $i'\le m$ if $\tau'$ selects the atom, and $i\le j\le m$; let $\xi\le\delta c$ be the box coordinate of the atom. We have the following dichotomy.
\begin{itemize}
\item[(a)] If $\tau'$ selects the atom with positive probability on $F$, then $Q(\xi,i)<Q(\xi,i')$ since $Q$ is strictly increasing in the rank for ranks up to $m$. The loss under $\bx$ is then smaller on a set of positive probability and not larger elsewhere, so $L_{\bx}(\tau')<L_{\bx'}(\tau')$, which is a contradiction.
\item[(b)] Otherwise $\tau'$ rejects the atom a.s.\ on $F$. We show that selecting the atom is strictly better under $\bx$ than continuing optimally, for small $\delta$. If $\bx\in\mathcal S$, the loss of selection is at most $Q(\delta c,m)\to Q(0,m)<Q_\infty$ as $\delta\to0$ ($Q(0,m)$ is $m$ in $\mathrm A_m$ and $0$ in $\mathrm B_m$). After the arrival all components of the state in box coordinates are at most $\delta c$, so by (i) and continuity of $v$ at $\bxero$ the continuation value is at least $v(\delta c,\dots,\delta c)\to Q_\infty$. If $\bx=\binf$, then $i=1$, the loss of selection is at most $Q(\delta c,1)\to Q(0,1)$, and the continuation value is at least $v(\binf)$ by (i). Here $v(\binf)>Q(0,1)$: by Theorem~\ref{thm:bellman}(iv) the value is attained by an optimal strategy, which loses at least $Q(\xi,1)>Q(0,1)$ when it selects an atom of level $\xi>0$, and $Q_\infty>Q(0,1)$ when it selects nothing. Modifying $\tau'$ to select the atom on $F$ (an event determined at the arrival of the atom) lowers the expected loss under $\bx$, so $\tau'$ is not optimal for $\bx$, again a contradiction.
\end{itemize}
Hence $v(\bx)<v(\bx')$.

(iii) Since $\bx\ge(x_1,\dots,x_1)$, part (i) and Theorem~\ref{thm:bellman}(iv) give $v(\binf)\le v(\bx)\le v(x_1,\dots,x_1)\to v(\binf)$ as $x_1\to\infty$.
\end{proof}

\subsection{The adaptive threshold}\label{sec:mon}
Suppose that at time $t$ the history in the original coordinates is $\by$. By Lemma~\ref{lem:self} and Theorem~\ref{thm:bellman}, an atom $(t,x)$ with initial rank $i=i(\by,x)$ is acceptable for the optimal strategy if and only if
\begin{equation}\label{eq:stopcond}
Q(\bar tx,\,i)\ \le\ v\big(\bar t\,(\by\cup x)\big).
\end{equation}
Let $H(t,\by)$ be the supremum of the scores $x$ satisfying (\ref{eq:stopcond}); this is a deterministic function of $t$ and $\by$.

We now define the random {\it adaptive threshold} process $(H_t,~t\in[0,1))$ on the whole unit interval. Starting from the initial history $\by_0$ (which may be empty), 
consider the reduced history updated only by the atoms arriving after time $0$ {\it above} the threshold: the arrivals below the threshold are acceptable, but they are ignored in the definition. Denoting 
$\by^\circ_t$ the history so obtained, we define 
$$H_t:=H(t,\by^\circ_t),~~~C_t:=v(\bar t\by^\circ_t).$$
Thus $C_t$ is the minimal expected loss to go from time $t$ for the history $\by^\circ_t$, whose entries are scores of the initial history or of atoms which were unacceptable at arrival, i.e.
the continuation value given that no acceptable atom has arrived; up to time $\tau$ this is the actual continuation value.
The process $(H_t)$ is  adapted to the arrivals of the PPP above the threshold itself, and the optimal strategy stops at the first atom below the threshold,
$$\tau=\inf\{t>0:\ \Pi\cap\{(t,x):x<H_t\}\ne\varnothing\}$$
($\tau=1$ if no such $t$).
Up to time $\tau$ the history $\by^\circ_t$ is the actual history. After $\tau$ the process $(H_t)$ is a formal continuation, which is convenient since it lives on the whole interval, and is the analogue of the memoryless rules, which stop at the first atom below a nonrandom precommitted threshold.

\begin{proposition}\label{prop:mon}
\begin{itemize}
\item[{\rm(i)}] The set of scores $x$ satisfying {\rm (\ref{eq:stopcond})} is an interval $[0,H_t)$ (possibly with the right endpoint).
\item[{\rm(ii)}] Along every path of the PPP the continuation value $C_t$ is nondecreasing in $t$; in particular $C_t\ge C_0=v(\by_0)$, which is $v(\binf)$ for the empty initial history.
\item[{\rm(iii)}] The threshold $H_t$ is nondecreasing in $t$: between arrivals it is continuous and either increases or is flat, and at an arrival of an atom above the threshold it either does not change or jumps up.
\item[{\rm(iv)}] If $H_t$ reaches a past score $y^\circ_i$, either continuously or by a jump at an arrival (then $y^\circ_i$ may be the score of the arriving atom), it stays flat at $y^\circ_i$ for a positive time, almost surely.
\end{itemize}
\end{proposition}
\begin{proof}
Write $\Phi_t(x):=Q(\bar tx,i)-v\big(\bar t(\by^\circ_t\cup x)\big)$ with $i=i(\by^\circ_t,x)$, so that $x$ is acceptable if and only if $\Phi_t(x)\le0$.

(i) Let $z<x$ and suppose that $x$ satisfies (\ref{eq:stopcond}). Then $i(\by^\circ_t,z)\le i(\by^\circ_t,x)$, $\bar tz<\bar tx$, and $Q$ is nondecreasing in $i$ and in the level; also $\by^\circ_t\cup z\le \by^\circ_t\cup x$ componentwise, so by Lemma~\ref{lem:mon} $v(\bar t(\by^\circ_t\cup z))\ge v(\bar t(\by^\circ_t\cup x))$, thus $z$ satisfies (\ref{eq:stopcond}) too.

(ii) Between arrivals above the threshold the reduced history $\by^\circ_t$ is fixed and $\bar t$ decreases, so $\bar t\by^\circ_t$ decreases componentwise and $C_t=v(\bar t\by^\circ_t)$ does not decrease by Lemma~\ref{lem:mon}. When an atom above the threshold arrives at time $t$, the new history $\by^\circ_t$ consists of the $m$ lowest scores among $\by^\circ_{t-}$ and the new score, so $\by^\circ_t\le\by^\circ_{t-}$ and $C_t$ jumps up (or does not change, if the new score is above the $m$-th lowest score of $\by^\circ_{t-}$).

(iii) {Between arrivals} the history $\by^\circ_t$ and the ranks $i(\by^\circ_t,x)$ are fixed, and $\Phi_t(x)$ is continuous and strictly decreasing in $t$, since $Q$ is increasing in the level $\bar tx$, which decreases, and $v$ increases when the arguments decrease. Hence the set of acceptable scores grows with $t$, which means that $H_t$ is nondecreasing; it moves continuously, since for the losses {\rm(\ref{QA})}, {\rm(\ref{QB})} the function $\Phi_t(x)$ is jointly continuous in $(t,x)$ and strictly increasing in $x$ between consecutive past scores, and jumps up at the past scores. 

{ At an arrival} of an atom at time $t$ with the score $x_0\ge H_{t-}$ (the arrivals below the threshold are ignored) the reduced history changes from $\by^\circ_{t-}$ to $\by^\circ_t$, the $m$ lowest among $\by^\circ_{t-}$ and $x_0$; this is no change if $x_0$ is above the $m$-th lowest score of $\by^\circ_{t-}$. Let $z<H_{t-}$. Since $x_0>z$, the initial rank of $z$ is the same under $\by^\circ_{t-}$ and $\by^\circ_t$, while $\by^\circ_t\cup z\le\by^\circ_{t-}\cup z$, so $v(\bar t(\by^\circ_t\cup z))\ge v(\bar t(\by^\circ_{t-}\cup z))$ and $\Phi(z)$ does not increase. Thus every $z<H_{t-}$ remains acceptable, so $H$ does not decrease. Here it is essential that the arrivals below the threshold are ignored, since such an arrival would raise the rank of the higher scores and lower $H_t$.

(iv) Suppose that $H_t=y^\circ_i$ at some time $t$, where the scores below $y^\circ_i$ are acceptable and those above are not. Compare the adapted loss at the levels $y^\circ_i-$ and $y^\circ_i+$: the initial rank is $i$ below and $i+1$ above, so the loss jumps up by $\Delta_i(\xi)>0$, which is $1$ in $\mathrm A_m$ and $e^{-\xi}\xi^{m-i}/(m-i)!$ in $\mathrm B_m$ for $i\le m$, whereas the continuation value is continuous in $x$. Thus the acceptable set ends at $y^\circ_i$ whenever $Q(\bar ty^\circ_i,i)\le v<Q(\bar ty^\circ_i,i+1)$, with $v$ the continuation value at $x=y^\circ_i$. If $H$ reaches $y^\circ_i$ continuously, then $Q(\bar ty^\circ_i,i)=v$ at the time of arrival, so the second inequality holds with the margin $\Delta_i$. If $H$ lands on $y^\circ_i$ by a jump, at the arrival of an atom $(t,x_0)$, the second inequality holds strictly, except in the borderline case $Q(\bar ty^\circ_i,i+1)=v$, in which $H$ would leave $y^\circ_i$ at once. The borderline case has probability zero: between consecutive arrivals, for a fixed score $x_0$ the difference $Q(\bar ty^\circ_i,i+1)-v$ is strictly decreasing in $t$, because $Q(\,\cdot\,,i+1)$ is strictly increasing in the level for $i+1\le m$ and $v$ does not decrease in $t$ by Lemma~\ref{lem:mon}(i) (for $i=m$ there is no borderline case, since $Q(\,\cdot\,,m+1)\ge Q_\infty>v$), so it vanishes at most once; the pairs $(t,x_0)$ of this kind form a Lebesgue null set, and by the compensator formula the expected number of atoms of $\Pi$ in it is zero. Between arrivals $Q(\bar ty^\circ_i,i+1)-v$ decreases continuously in $t$, hence the inequality persists, and $H$ stays flat at $y^\circ_i$, during a positive time interval.
\end{proof}

Thus if the threshold is not in a flat stretch, it drifts smoothly between arrivals in a manner determined by the frozen history, and jumps up when an atom arrives above it and lowers the history. A jump takes the threshold up to a level where the loss of acceptance and the continuation value balance, or to a past score, and in the latter case the threshold stays there by (iv). Since the loss above a past score is higher by a jump, we do not expect the threshold to be carried across the score of the arriving atom; 
for ${\rm A}_m$ this is Lemma~\ref{lem:across} below, which rests on the linearity of the loss, below truncation, in the rank;
for ${\rm B}_2$ this is verified in Section~\ref{sec:two}.
For $m=1$ the drift follows a hyperbola up to the lowest score and is frozen afterwards, but for $m\ge2$ the stopping boundary in box coordinates depends on the state, which itself shrinks in time between history updates. Every time the threshold reaches a past score
(either by drift or a jump) it stays flat there, as the loss jumps above the atom. In words, (ii) says that {\it if no atom has been accepted up to time $t$, then the prospects of the observer can only have worsened}.

\begin{lemma}\label{lem:across}
In $\mathrm A_m$, for every state $\bx\in\mathcal S$ and every $\xi<x_m$,
$$v(\bx\cup\xi)\ \le\ v(\bx)+1 .$$
Consequently, when a low atom $(t,x_0)$ which is not acceptable arrives, no score above $x_0$ is acceptable after the arrival: the threshold is not carried across $x_0$.
\end{lemma}
\begin{proof}
Let $\tau$ be a strategy in the problem with the history $\bx$, and apply it with the history $\bx\cup\xi$, as in the proof of Lemma~\ref{lem:mon}(i). Compared with the levels $x_1,\dots,x_m$, the $m+1$ levels $x_1,\dots,x_m,\xi$ raise the initial rank $i$ of every atom by at most one; removing the largest of them, $x_m$, changes only the ranks which exceed $m$, so it does not change $\min(i,m+1)$. Hence the loss {\rm(\ref{QA})} of the selected atom is larger by at most $1$ under $\bx\cup\xi$, the loss of no selection is the same, and $L_{\bx\cup\xi}(\tau)\le L_{\bx}(\tau)+1$. Taking the infimum over $\tau$ gives the inequality.

Let $\by$ be the history just before the arrival of $(t,x_0)$, let $i_0=i(\by,x_0)$ be the initial rank of the atom and $\by'=\by\cup x_0$ the history after the arrival. The atom is not acceptable, that is, $Q(\bar tx_0,i_0)>v(\bar t\by')$. A score $z\ge y'_m$ is not low, hence not acceptable. For $x_0<z<y'_m$ the initial rank of $z$ under $\by'$ is at least $i_0+1\le m$, so by the monotonicity of $Q$ and the first part of the proof
$$Q\big(\bar tz,i(\by',z)\big)\ \ge\ Q(\bar tx_0,i_0)+1\ >\ v(\bar t\by')+1\ \ge\ v\big(\bar t(\by'\cup z)\big),$$
and $z$ does not satisfy {\rm(\ref{eq:stopcond})}.
\end{proof}

\begin{remark}[Threshold monotonicity in discrete time]\label{rem:BF}
For the Robbins problem Bruss and Ferguson \cite{BFhp} (Lemma~4.1) assert that the optimal stopping time is a stepwise monotone increasing (adaptive) threshold strategy. In the proof (p.~13) they use the almost sure
{\it pathwise} continuation value monotonicity, which, however,  contradicts the stopping value
monotonicity in $n$ (indeed, the first `extremely risky' arrival is rejected and the problem is reduced to more favourable case $n-1$); thus the threshold monotonicity remains a conjecture. 
 A direct computation for $n\le7$ finds the inequality (4.10) of \cite{BFhp} violated in about $1.4\%$ of the states examined, whereas the conclusion of the lemma was verified throughout.
This nuisance does not appear in the PPP problem, since
 the continuation value increases pathwise along every path of the reduced history by the scaling property (Proposition~\ref{prop:mon}(ii)); the monotonicity of the threshold process was noted in \cite{GnedinJAP} (Lemma~6).
\end{remark}

\begin{remark}[Rank-information strategies]
The  $i$-records arrive at times of Poisson process with intensity $1/t$, independently across $i$
(this is a version of seminal Ignatov's Theorem, see \cite{GnedinCorners}). 
The superposition of these point processes is equivalent to  the infinite-arrivals model of \cite{GiSamuels}, thus
the strategy optimal in the class of  rules adapted to only initial ranks 
is the same as in  \cite{GiSamuels, GnedinUniform}. 
This strategy is defined by 
 cutoffs $0\leq t_1\leq\dots\leq t_m\le t_{m+1}=1$, and accepts the first atom with initial rank at most $i$ arriving within a slot  $[t_i,t_{i+1})$. 
An atom arriving at time $s\in[t_i,t_{i+1})$ has initial rank at most $i$ if and only if its score is below $y_i(s-)$, the $i$-th lowest score observed before $s$.  
Arrivals above $y_i$ do not change $y_i$,  so the threshold  is $0$ before $t_1$, is flat at the level $y_i(t_i)$ on $[t_i,t_{i+1})$, and jumps up at the cutoff times $t_{i+1}$ to  $y_{i+1}> y_i$
until an acceptable arrival. Ignoring acceptable arrivals yields a nondecreasing adaptive threshold process.
To compare, the optimal `full-information' threshold process $(H_t)$ may drift and jump on any time interval.

\end{remark}

\section{One tracked score: $m=1$}\label{sec:one}
For $m=1$ only records can be accepted, and the state is a single number: the box coordinate $\xi$ of the lowest score $Y$ observed so far. In the original coordinates we call $(t,1]\times(0,Y)$ the {\it south-east window} of the current record. The area is $\xi=\bar tY$; the new records can only arrive into this window, other arrivals being of no interest. The area decreases with time, and by the first arrival into the window it is replaced by the smaller window below the record.

\subsection{A hyperbolic threshold}
Suppose a record arrives with box coordinate $\xi$. If it is accepted, the loss is $Q(\xi,1)$, which is increasing in $\xi$. If it is continued, the state becomes $\xi$, and the future loss is the continuation value $v(\xi)$, which is decreasing in $\xi$ by Lemma~\ref{lem:mon}. Therefore
$$\{\xi:\ Q(\xi,1)\le v(\xi)\}=(0,b]$$
for some $b>0$: {\it the optimal strategy stops at the first record $(t,x)$ with $\bar tx\le b$, that is, at the first record below the hyperbola $h(t)=b/\bar t$.}
We call this strategy $\mathrm R_b$, and the problem is reduced to optimisation over $b$. In the original PPP coordinates the adaptive threshold is $H_t=\min(Y,\, b/\bar t)$, where $Y$ is the lowest score of the reduced history (up to the time $\tau$, the lowest score observed): it drifts along the hyperbola until it meets $Y$, where it stays frozen.

\begin{figure}[ht]
\centering
\begin{tikzpicture}
\begin{axis}[width=0.78\textwidth,height=0.46\textwidth,xmin=0,xmax=0.9,ymin=0,ymax=4,
 xlabel={time $t$},ylabel={score $x$},ytick={0,1,2,3,4},xtick={0,0.1,0.2,0.3,0.4,0.5,0.6,0.7,0.8,0.9},
 xtick pos=bottom,ytick pos=left,tick label style={font=\footnotesize},label style={font=\small},axis line style={black},clip=true,
 scaled ticks=false]
\fill[gray!22] (axis cs:0,0) rectangle (axis cs:0.5566,1.3930);
\addplot[black,thick] coordinates {(0.0000,0.6176) (0.0058,0.6212) (0.0116,0.6249) (0.0174,0.6286) (0.0231,0.6323) (0.0289,0.6360) (0.0347,0.6399) (0.0405,0.6437) (0.0463,0.6476) (0.0521,0.6516) (0.0579,0.6556) (0.0636,0.6596) (0.0694,0.6637) (0.0752,0.6679) (0.0810,0.6721) (0.0868,0.6763) (0.0926,0.6806) (0.0983,0.6850) (0.1041,0.6894) (0.1099,0.6939) (0.1157,0.6985) (0.1215,0.7031) (0.1273,0.7077) (0.1331,0.7124) (0.1388,0.7172) (0.1446,0.7221) (0.1504,0.7270) (0.1562,0.7320) (0.1620,0.7370) (0.1678,0.7422) (0.1736,0.7474) (0.1793,0.7526) (0.1851,0.7580) (0.1909,0.7634) (0.1967,0.7689) (0.2025,0.7745) (0.2083,0.7801) (0.2141,0.7859) (0.2198,0.7917) (0.2256,0.7976) (0.2314,0.8036) (0.2372,0.8097) (0.2430,0.8159) (0.2488,0.8222) (0.2546,0.8286) (0.2603,0.8350) (0.2661,0.8416) (0.2719,0.8483) (0.2777,0.8551) (0.2835,0.8620) (0.2893,0.8690) (0.2950,0.8761) (0.3008,0.8834) (0.3066,0.8908) (0.3124,0.8983) (0.3182,0.9059) (0.3240,0.9136) (0.3298,0.9215) (0.3355,0.9295) (0.3413,0.9377) (0.3471,0.9460) (0.3529,0.9545) (0.3587,0.9631) (0.3645,0.9719) (0.3703,0.9808) (0.3760,0.9899) (0.3818,0.9991) (0.3876,1.0086) (0.3934,1.0182) (0.3992,1.0280) (0.4050,1.0380) (0.4108,1.0482) (0.4165,1.0586) (0.4223,1.0692) (0.4281,1.0800) (0.4339,1.0910) (0.4397,1.1023) (0.4455,1.1138) (0.4512,1.1255) (0.4570,1.1375) (0.4628,1.1498) (0.4686,1.1623) (0.4744,1.1751) (0.4802,1.1882) (0.4860,1.2015) (0.4917,1.2152) (0.4975,1.2292) (0.5033,1.2435) (0.5091,1.2582) (0.5149,1.2732) (0.5207,1.2886) (0.5265,1.3043) (0.5322,1.3204) (0.5380,1.3370) (0.5438,1.3539) (0.5496,1.3713) (0.5554,1.3892) (0.5612,1.4075) (0.5670,1.4263) (0.5727,1.4456) (0.5785,1.4654) (0.5843,1.4858) (0.5901,1.5068) (0.5959,1.5284) (0.6017,1.5506) (0.6074,1.5734) (0.6132,1.5969) (0.6190,1.6212) (0.6248,1.6462) (0.6306,1.6720) (0.6364,1.6986) (0.6422,1.7260) (0.6479,1.7544) (0.6537,1.7837) (0.6595,1.8140) (0.6653,1.8454) (0.6711,1.8778) (0.6769,1.9115) (0.6827,1.9463) (0.6884,1.9824) (0.6942,2.0199) (0.7000,2.0589) (0.7058,2.0994) (0.7116,2.1415) (0.7174,2.1853) (0.7232,2.2310) (0.7289,2.2786) (0.7347,2.3283) (0.7405,2.3802) (0.7463,2.4345) (0.7521,2.4913) (0.7579,2.5508) (0.7637,2.6133) (0.7694,2.6788) (0.7752,2.7478) (0.7810,2.8204) (0.7868,2.8969) (0.7926,2.9777) (0.7984,3.0631) (0.8041,3.1536) (0.8099,3.2496) (0.8157,3.3516) (0.8215,3.4602) (0.8273,3.5762) (0.8331,3.7001) (0.8389,3.8329) (0.8446,3.9757)};
\addplot[only marks,mark=*,mark size=1.4pt,black] coordinates {(0.2227,2.2860) (0.4203,1.2186) (0.6097,2.7638) (0.6762,1.5656) (0.8138,0.1033)};
\addplot[only marks,mark=*,mark size=1.4pt,black] coordinates {(0.5566,0.6012)};
\draw[black,thick] (axis cs:0.5566,0.6012) circle (3.2pt);
\node[font=\small,anchor=west] at (axis cs:0.03,3.55) {curve $x=b/\bar t$};
\end{axis}
\end{tikzpicture}
\caption{A realisation of $\Pi$ with the hyperbola $x=b/\bar t$, $b=0.6176$, the optimal threshold in problem $\mathrm A_1$. The atom at $(0.42,1.22)$ arrives above the curve; the curve passes its level at $t=0.49$, so the window meets this atom through its upper side. The first record below the curve afterwards (circled, $t=0.557$) is selected. Shaded: the south-west window with the apex on the hyperbola at the time of selection; the selected atom is on its right side.}\label{fig:hyp}
\end{figure}
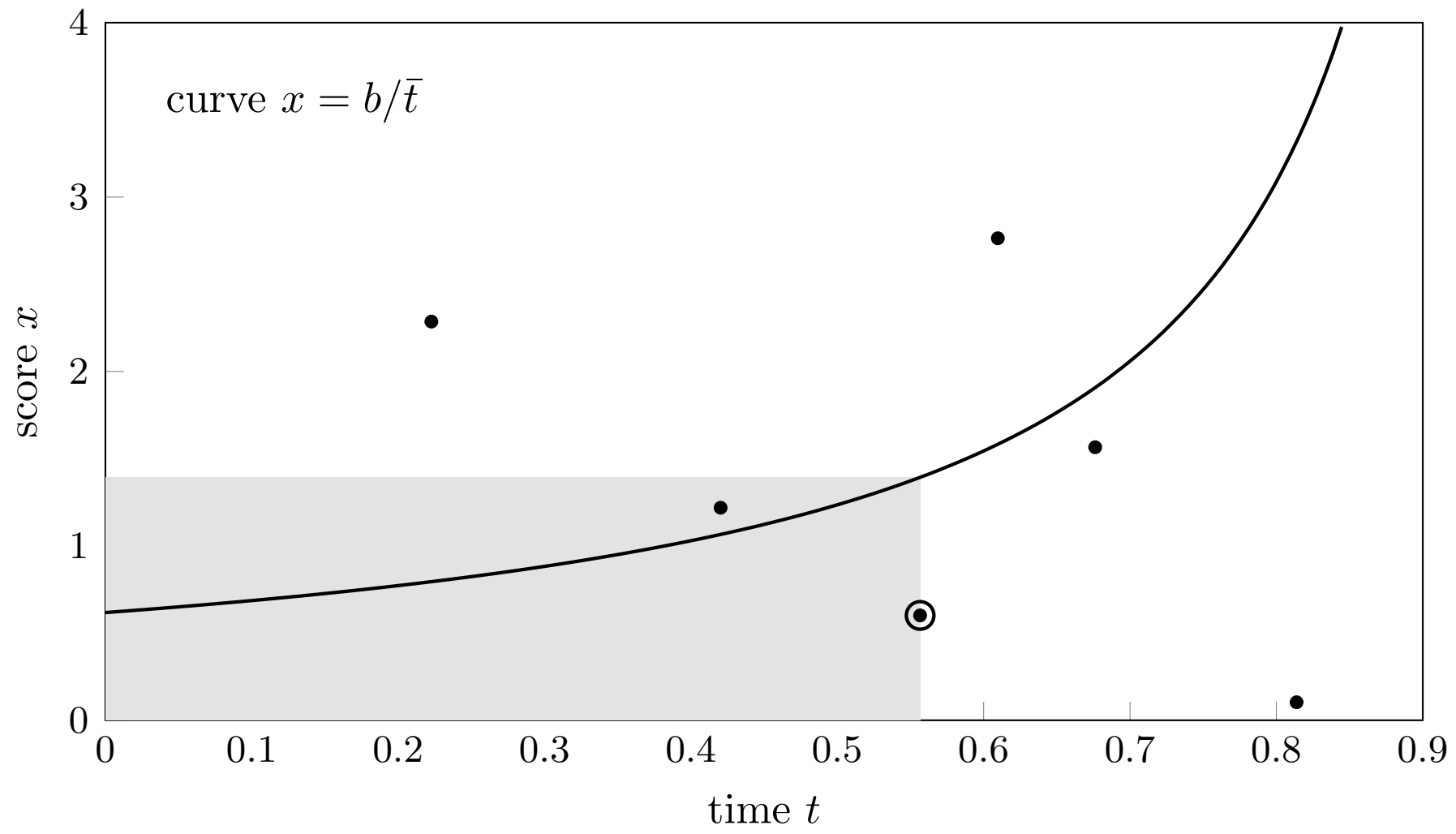

\subsection{The south-west window and the loss of $\mathrm R_b$}
Let $h(t)=b/\bar t$. Consider the { south-west window}
$$W_t=(0,t]\times(0,h(t)),$$
with the apex $(t,h(t))$ on the hyperbola, and let the apex slide along the hyperbola, so that its time coordinate moves at unit speed. The window grows, its area is $|W_t|=t\,h(t)=bt/\bar t$, and the new area is swept by two sides: the right side, of length $h(t)$, and the upper side, of length $t$, rising at the speed $h'(t)$. The atoms of $\Pi$ are met by the window according to a Poisson process in $t$, with the intensity $h(t)+t h'(t)=\dd|W_t|/\dd t$, which splits in the rates $h(t)$ for the right side and $t h'(t)$ for the upper side. The probability that no atom is met before time $t$ is $e^{-bt/\bar t}$. To assess the rule we need to consider two cases.

\begin{itemize}
\item[\rm (right)] {\it The first atom is met by the right side} at time $t$: this is an atom $(t,x)$ with $x$ uniform on $(0,h(t))$. Since the window was empty at time $0$, this is a record below the hyperbola, so the strategy $\mathrm R_b$ stops with the box coordinate $\bar tx$  uniform on $(0,b)$, and the expected loss is
$$\bar g(b):=\frac1b\int_0^bQ(u,1)\,\dd u .$$
\item[\rm (upper)] {\it The first atom is met by the upper side} at time $t$: this is an earlier atom $(s,h(t))$, $s<t$, located above the hyperbola at the time $s$ of its arrival and was not accepted. The scores observed before $t$  are all above $h(t)$, so it is the lowest  so far, and its south-east window has the area $\bar th(t)=b$. From time $t$ the strategy selects greedily the next record (if any),
 as  records fall into the window below the hyperbola. The expected loss of this future action is $G(b)$ computed below.
\end{itemize}

The probabilities of the two cases are
$$
\P[\text{right}]=\int_0^1 e^{-bt/\bar t}\,\frac{b}{\bar t}\,\dd t=be^bJ_1(b),\qquad \P[\text{upper}]=\int_0^1 e^{-bt/\bar t}\,\frac{bt}{\bar t^{\,2}}\,\dd t=1-be^bJ_1(b),
$$
where $J_1(b):=\int_b^\infty e^{-s}s^{-1}\dd s$ is the exponential integral function.

\paragraph{The greedy phase.} Suppose the state is $\xi$, that is, at time $0$ the lowest score is $\xi$ and the south-east window $(0,1]\times(0,\xi)$ has the area $\xi$.
With probability $e^{-\xi}$  the window contains no atoms. Otherwise the first atom arrives at  time $s\in(0,1)$ with exponential  density $\xi e^{-\xi s}$, it is then a record with the score $\xi u$, where $u$ is sampled independently from  ${\rm U}[0,1]$, and the new state becomes $\xi(1-s)u$.
It follows that
\begin{equation}\label{eq:G}
G(\xi)=Q_\infty e^{-\xi}+\int_0^\xi e^{y-\xi}\,\bar g(y)\,\dd y,\qquad\text{equivalently}\qquad G'+G=\bar g,\quad G(0)=Q_\infty.
\end{equation}
For $\mathrm A_1$ and $\mathrm B_1$, with $J_2(\xi):=\int_0^\xi (e^s-1)s^{-1}\dd s$,
\begin{align}
&\mathrm A_1:&& g(\xi)=1+\xi,\ Q_\infty=2: && G(\xi)=\frac{1+\xi}{2}+\frac32e^{-\xi},\label{GA}\\
&\mathrm B_1:&& g(\xi)=1-e^{-\xi},\ Q_\infty=1: && G(\xi)=1-e^{-\xi}J_2(\xi). \label{GB}
\end{align}
In detail, for $\mathrm A_1$, $\bar g(y)=1+y/2$ and $\int_0^\xi e^y(1+ y/2)\dd y=e^\xi-1+\frac12((\xi-1)e^\xi+1)$; and for $\mathrm B_1$, $\bar g(y)=1-(1-e^{-y}){y}^{-1}$ and $\int_0^\xi e^y (1-e^{-y}){y}^{-1}\dd y=J_2(\xi)$.

\begin{proposition}\label{prop:Rb}
The expected loss of the rule $\mathrm R_b$ started from $\binf$ is
$$R(b)=be^bJ_1(b)\,\bar g(b)+\big(1-be^bJ_1(b)\big)\,G(b)=G(b)+b\,e^{b}J_1(b)\,G'(b).$$
\end{proposition}
\begin{proof}
The first expression is the decomposition into the cases of the right and the upper side. The second follows from $\bar g=G+G'$.
\end{proof}

For $\mathrm B_1$ the two terms are the contributions of the upper and the right crossing of the window to the success probability, which is the splitting used for the best-choice problem in \cite{Memoryless, GKS}; Proposition~\ref{prop:Rb} extends it to an arbitrary loss $Q(\xi,1)$.

\subsection{Optimising $b$}
By (\ref{eq:G}), $\xi(G'+G)=\int_0^\xi g(z){\rm d}z$, and differentiating, $\xi G''+(1+\xi)G'=g-G$. With $\rho(b)=be^bJ_1(b)$ we have $\rho'(b)=(1+b)e^bJ_1(b)-1$, so
\begin{equation}\label{eq:dR}
R'(b)=G'(1+\rho')+\rho G''=e^bJ_1(b)\big[(1+b)G'(b)+bG''(b)\big]=e^{b}J_1(b)\,\big(Q(b,1)-G(b)\big).
\end{equation}
The right side vanishes at the {\it balance} point: accepting a record of level $b$ and waiting for the next record are equally good. For $\mathrm A_1$, $Q(b,1)-G(b)=\frac{1+b}{2}-\frac32e^{-b}$ is increasing; for $\mathrm B_1$, $Q(b,1)-G(b)=e^{-b}(J_2(b)-1)$ has the sign of $J_2(b)-1$, and $J_2$ is increasing. In both cases the difference is negative at $b=0$ and has a single zero, so $R$ decreases up to this root and increases afterwards.

\begin{theorem}\label{thm:one}
{\rm (i)} In ${\mathrm A}_1$ the optimal strategy stops at the first record below the hyperbola $x=b_{\rm A}/\bar t$, where $b_{\rm A}=0.617642\ldots$ is the root of $(1+b)e^{b}=3$. The value is
$${\bf a}_1=1+b_{\rm A}-K_{\rm A}J_1(b_{\rm A})=1.462488\ldots,\qquad K_{\rm A}=\tfrac12b_{\rm A}^2e^{b_{\rm A}}=0.353739\ldots .$$
{\rm (ii)} In ${\mathrm B}_1$ the optimal strategy stops at the first record below the hyperbola $x=b_{\rm B}/\bar t$, where $b_{\rm B}=0.804352\ldots$ is the root of $J_2(b)=1$. The value is
$${\bf b}_1=1-e^{-b_{\rm B}}-K_{\rm B}J_1(b_{\rm B})=0.419836\ldots,\qquad K_{\rm B}=e^{b_{\rm B}}-1-b_{\rm B}=0.430896\ldots,$$
so the maximal probability of success is $1-{\bf b}_1=0.580164\ldots$.
\end{theorem}
\begin{proof}
The optimal strategy is $\mathrm R_b$ for some $b$, and by (\ref{eq:dR}) the best $b$ is the root of the balance equation $Q(b,1)=G(b)$. For $\mathrm A_1$ this is $\frac{1+b}2=\frac32e^{-b}$, that is $(1+b)e^b=3$; for $\mathrm B_1$ it is $J_2(b)=1$. At the root $G(b)=Q(b,1)$ equals $1+b$, respectively $1-e^{-b}$, and $-be^bG'(b)$ equals $\frac12b^2e^b$ (using $e^b=3/(1+b)$), respectively $e^b-1-b$, since $G'=e^{-b}(1-\frac{e^b-1}{b})$ when $J_2(b)=1$.
\end{proof}
The second case is the known solution of the best-choice problem \cite{GM, BC-SPA, GKS}. We have used the optimality  equation only through the interval structure of the optimal strategy.

\begin{remark}[The value function]\label{rem:vone}
For $m=1$ equation (\ref{eq:bellman}) has the form 
$$\xi v'(\xi)+\xi v(\xi)=\int_0^\xi\min\{Q(w,1),v(w)\}\dd w.$$
For $\xi\le b$ it gives $v=G$, and for $\xi>b$ the integrand is $v$, so $\xi v''+(1+\xi)v'=0$, that is $(\xi e^{\xi}v')'=0$, and $v'$ is continuous at $b$. Hence
$$v(\xi)=G(\xi),\ \ \xi\le b;\qquad v(\xi)=v(\binf)+K_*J_1(\xi),\ \ \xi\ge b,$$
with $K_*=K_{\rm A}$ for $\mathrm A_1$ and $K_*=K_{\rm B}$ for $\mathrm B_1$. Here $v(\binf)={\bf a}_1$, respectively ${\bf b}_1$. This function is the edge function of the two-dimensional problems below.
\end{remark}

\section{Two tracked scores: $m=2$}\label{sec:two}
For $m=2$ the state $\bx=(x_1,x_2)$ is the lowest and the second lowest scores in box coordinates. An arriving atom at level $\xi$ is a record if $\xi<x_1$ (the new state is $(\xi,x_1)$) and a $2$-record if $x_1<\xi<x_2$ (the new state is $(x_1,\xi)$); atoms with $\xi>x_2$ are ignored. Equation (\ref{eq:bellman}) takes the form
\begin{equation}\label{eq:two}
x_1\frac{\partial v}{\partial x_1}+x_2\frac{\partial v}{\partial x_2}+x_2\,v=\int_0^{x_1}\min\{Q(\xi,1),\,v(\xi,x_1)\}\,\dd\xi+\int_{x_1}^{x_2}\min\{Q(\xi,2),\,v(x_1,\xi)\}\,\dd\xi .
\end{equation}

\subsection{Two boundaries}
In the state $(x_1,x_2)$ a record of level $\xi$ is acceptable by the optimal strategy if and only if $Q(\xi,1)\le v(\xi,x_1)$, and a $2$-record of level $\xi$ is acceptable if and only if $Q(\xi,2)\le v(x_1,\xi)$. By Proposition~\ref{prop:mon}(i) the acceptable levels form an interval, and both boundaries of the stopping sets are determined by {\it balancing} the loss of acceptance against the continuation value. Therefore, for a given $x_1$ there are two cases.
\begin{itemize}
\item[{\rm(i)}] {\it $2$-records are acceptable.} Then also records are acceptable, and the $2$-records are acceptable up to the level $\beta_2(x_1)>x_1$, the {\it rank-$2$ boundary}, which is the root of
\begin{equation}\label{eq:mu}
Q(\beta_2,2)=v(x_1,\beta_2).
\end{equation}
\item[{\rm(ii)}] {\it $2$-records are not acceptable}. Then the records are acceptable up to the level $\beta_1(x_1)\le x_1$, the {\it rank-$1$ boundary}, which is the root of
\begin{equation}\label{eq:rho}
Q(\beta_1,1)=v(\beta_1,x_1),
\end{equation}
or $\beta_1(x_1)=x_1$ if all records are acceptable.
\end{itemize}
Both equations are {\it balance} conditions (continuous fit): the left side is the loss of acceptance, the right side is the loss of continuation after the new state is formed. In the state space the sets $\{x_1\le\beta_1(x_2)\}$ and $\{x_2\le\beta_2(x_1)\}$ consist of the states which are reached by an acceptable record, respectively $2$-record (Figure~\ref{fig:stop}). In the original coordinates, with $y^\circ_1<y^\circ_2$ the two lowest scores of the reduced history, the adaptive threshold is
$$
H_t=\min\big(y^\circ_2,\ \beta_2(\bar ty^\circ_1)/\bar t\big)\ \text{ in the case (i)},\qquad H_t=\beta_1(\bar ty^\circ_1)/\bar t\ \text{ in the case (ii)}.
$$
Between arrivals $H_t$ is nondecreasing by Proposition~\ref{prop:mon}, but it is not a hyperbola, because $\beta_1,\beta_2$ are evaluated at the shrinking state $\bar ty^\circ_1$, while the past scores are frozen. The threshold  stays flat at $y^\circ_2$, respectively $y^\circ_1$, when it reaches the past score.

\begin{remark}[Jumps for $m=2$]\label{rem:jump2}
In the case (i) every record is below $y^\circ_1\le H_t$ and is accepted, while a $2$-record above $H_t$ can only lower $y^\circ_2$ to its score $x_0\ge H_t$, which leaves $H_t=\min(y^\circ_2,\beta_2(\bar ty^\circ_1)/\bar t)$ unchanged. In the case (ii) an atom above $y^\circ_1$ does not change $H_t$, and a rejected record $x_0\in[H_t,y^\circ_1)$ becomes the new lowest score (this requires $\bar ty^\circ_1>x^*$, since otherwise $H_t=y^\circ_1$). Then $\bar tx_0\ge\beta_1(\bar ty^\circ_1)\ge H_0>\gamma$. The first inequality holds because $x_0\ge H_t$. The second holds because $\beta_1$ is nonincreasing on $[x^*,\infty)$ with $\beta_1(\infty)=H_0$: the root of $Q(\xi,1)=v(\xi,u)$ moves to the left as $u$ grows, since $v$ decreases in its second argument (Lemma~\ref{lem:mon}). The third is Lemma~\ref{lem:b1gamma} (numerically $H_0=0.902$ and $\gamma=0.389$ in $\mathrm A_2$, $H_0=1.088$ and $\gamma=0.706$ in $\mathrm B_2$). So we stay in the case (ii), and the new threshold is $\beta_1(\bar tx_0)/\bar t$. It is at least $H_t$, because $\bar tx_0\ge\beta_1(\bar ty^\circ_1)$ while $\beta_1(u)=u$ for $u\le x^*$ and $\beta_1$ is nonincreasing on $[x^*,\infty)$, and it is at most $x_0$, because $\beta_1(u)\le u$. Thus a jump takes $H_t$ at most to the score of the arriving atom, and exactly there if $\bar tx_0\le x^*$ (all records are acceptable); then $H_t$ stays flat at $x_0$ by Proposition~\ref{prop:mon}(iv).
\end{remark}

\begin{lemma}\label{lem:shape}
In both problems $\mathrm A_2$ and $\mathrm B_2$ the function $x\mapsto Q(x,2)-v(x,x)$ is continuous and strictly increasing, with a unique zero $\gamma$. For $x_1<\gamma$ equation {\rm (\ref{eq:mu})} has a unique root $\beta_2(x_1)>x_1$, which is a strictly decreasing function of $x_1$ with $\beta_2(x_1)\to\gamma$ as $x_1\uparrow\gamma$; a $2$-record is acceptable if and only if its level is below $\beta_2(x_1)$. For $x_1\ge\gamma$ no $2$-record is acceptable.
\end{lemma}
\begin{proof}
$Q(x,2)$ is strictly increasing and $v(x,x)$ is strictly decreasing by Lemma~\ref{lem:mon}, and $v$ is continuous. As $x\to0$ the difference tends to $Q(0,2)-Q_\infty<0$ (it equals $2-3=-1$ in $\mathrm A_2$ and $0-1=-1$ in $\mathrm B_2$), and for large $x$ it is positive: $Q(x,2)=2+x\to\infty$ in $\mathrm A_2$, and $Q(x,2)=1-e^{-x}\to1>{\bf b}_2$ in $\mathrm B_2$. For fixed $x_1$ the function $\xi\mapsto Q(\xi,2)-v(x_1,\xi)$ is strictly increasing on $(x_1,\infty)$ and positive for large $\xi$; it is negative at $\xi=x_1$ if and only if $x_1<\gamma$. This gives the root and the threshold form. Raising $x_1$ lowers $v(x_1,\xi)$ for every $\xi$, so the root moves to the left.
\end{proof}
Thus the rank-2 boundary is a decreasing curve $x_2=\beta_2(x_1)$ in the plane, which runs from a point $(0,\beta_2(0))$ on the axis to the point $(\gamma,\gamma)$ on the diagonal.

\subsection{The edge}
The score $x_1\to0$ approaches  the lowest possible, which contributes $1$ to every other initial rank. 
The loss in $\mathrm A_2$ is $\min(i,3)+\xi=1+\min(i-1,2)+\xi$, thus the problem with the state $(0,x)$ amounts to problem $\mathrm A_1$ with the state $x$, and the loss increased by  $1$. In $\mathrm B_2$ the loss $\mathbf 1(r>2)$ becomes $\mathbf 1(r-1>1)$. Hence
$$
v_{\mathrm A_2}(0,x)=1+v_{\mathrm A_1}(x),\qquad v_{\mathrm B_2}(0,x)=v_{\mathrm B_1}(x).
$$
\begin{proposition}\label{prop:edge}
At the edge $x_1=0$ the rank-$2$ boundary is $\beta_2(0)=b_{\rm A}=0.617642\ldots$ in $\mathrm A_2$ and $\beta_2(0)=b_{\rm B}=0.804352\ldots$ in $\mathrm B_2$, which are the hyperbola constants of Theorem~{\rm \ref{thm:one}}.
\end{proposition}
\begin{proof}
For $\mathrm A_2$, $Q(\beta,2)=2+\beta$ and $v(0,\beta)=1+v_{\mathrm A_1}(\beta)$, so (\ref{eq:mu}) reads $1+\beta=v_{\mathrm A_1}(\beta)$, which is the balance equation $Q(\beta,1)=v_{\mathrm A_1}(\beta)$ defining $b_{\rm A}$. For $\mathrm B_2$, $Q(\beta,2)=1-e^{-\beta}$ is the loss of a record in $\mathrm B_1$, and $v(0,\beta)=v_{\mathrm B_1}(\beta)$, which defines $b_{\rm B}$.
\end{proof}
On the edge the stopping set for $2$-records is the  interval $[0,\beta_2(0)]$, and the value $v(0,x_2)$ is explicit for all $x_2$ by Remark~\ref{rem:vone}.

\subsection{The lower bound and the region $\mathcal R$}
Let $q(i):=Q(0,i)$, which is $\min(i,m+1)$ in $\mathrm A_m$ and $\mathbf 1(i>m)$ in $\mathrm B_m$; then $Q(\xi,i)\ge q(i)$ for all $\xi$ and $q(i)=Q_\infty$ for $i>m$. The following bound is the analogue, for the adapted losses, of the bound (14) in \cite{GnedinJAP} for the loss depending on the final rank, and is proved in the same way.
\begin{lemma}\label{lem:lb}
For every state $\bx$, with $x_0:=0$,
$$v(\bx)\ \ge\ \sum_{j=1}^{m}q(j)\big(e^{-x_{j-1}}-e^{-x_j}\big)+Q_\infty e^{-x_m}.$$
In particular, for $m=2$, $v(x_1,x_2)\ge q(1)(1-e^{-x_1})+q(2)(e^{-x_1}-e^{-x_2})+Q_\infty e^{-x_2}$.
\end{lemma}
\begin{proof}
In the problem with the initial history $\bx$, let $E_j$ be the event that $\Pi$ has no atom with score below $x_j$, so $\P[E_j]=e^{-x_j}$. On $E_j$ every accepted atom has the initial rank at least $j+1$. Hence the initial rank $I$ of the accepted atom (with $I=m+1$ if there is none) satisfies $\P[I\ge j+1]\ge e^{-x_j}$ for $j\le m$, and the loss is at least $q(I)$. By summation by parts, $\E q(I)=q(1)+\sum_{j\ge1}(q(j+1)-q(j))\P[I\ge j+1]$ is bounded below by the right side, since $q$ is nondecreasing.
\end{proof}

\begin{lemma}\label{lem:b1gamma}
In $\mathrm A_2$ and $\mathrm B_2$ we have $Q(\gamma,1)<v(\gamma,u)$ for every $u>\gamma$, with a difference bounded away from zero. Consequently $\beta_1(u)>\gamma$ for $u>x^*$, and $H_0>\gamma$.
\end{lemma}
\begin{proof}
In $\mathrm A_2$ Lemma~\ref{lem:lb} gives $v(\gamma,u)\ge1+e^{-\gamma}+e^{-u}>1+e^{-\gamma}$, while $Q(\gamma,1)=1+\gamma$, and $e^{-\gamma}=0.678>\gamma=0.389$.

In $\mathrm B_2$ consider the problem with the initial history $(\gamma,u)$, and let $E$ be the event that $\Pi$ has no atom with score below $\gamma$, so $\P[E]=e^{-\gamma}$. On $E$ the score $\gamma$ of the history is below every atom, hence a selected atom has final rank at most $2$ only if it is the lowest atom of $\Pi$. The event $E$ depends only on the atoms with score below $\gamma$, which are independent of the other atoms; so given $E$ the process $\Pi$ is a Poisson process of unit intensity on $[0,1]\times(\gamma,\infty)$, and after the shift of the scores by $\gamma$ the task is to select the lowest atom, which is problem $\mathrm B_1$ with empty history; every strategy fails in it with probability at least ${\bf b}_1$. Thus $v(\gamma,u)\ge e^{-\gamma}\,{\bf b}_1=0.207$ by Theorem~\ref{thm:one}, while $Q(\gamma,1)=1-(1+\gamma)e^{-\gamma}=0.158$.

Both bounds do not depend on $u$, so they hold for $v(\gamma,\infty)$ as well. Since $\xi\mapsto Q(\xi,1)-v(\xi,u)$ is increasing, its root $\beta_1(u)$ is above $\gamma$ for $u>x^*$, and so is the root $H_0$ of $Q(\xi,1)=v(\xi,\infty)$.
\end{proof}

Let $\mathcal R$ be the set of states $(x_1,x_2)$ with $x_1<\gamma$ and $x_1<x_2\le\beta_2(x_1)$, that is, the region between the diagonal and the rank-2 boundary. In $\mathcal R$ every $2$-record is acceptable. If also every record is acceptable, then (\ref{eq:two}) is linear:
\begin{equation}\label{eq:lin}
x_1\frac{\partial v}{\partial x_1}+x_2\frac{\partial v}{\partial x_2}+x_2\,v=\int_0^{x_1}Q(\xi,1)\,\dd\xi+\int_{x_1}^{x_2}Q(\xi,2)\,\dd\xi .
\end{equation}
The region $\mathcal R$ is star-shaped with respect to the origin, because $\beta_2$ is decreasing: if $(x_1,x_2)\in\mathcal R$ and $0<s<1$, then $sx_2<x_2\le\beta_2(x_1)<\beta_2(sx_1)$. Along a ray $x_2=\kappa x_1$ the left side of (\ref{eq:lin}) is an ordinary differential operator for $f(s)=v(s,\kappa s)$, namely $sf'+\kappa sf$, with the initial value $f(0)=Q_\infty$. This determines $v$ on $\mathcal R$. All records are acceptable in the state $(x_1,x_2)$ if $Q(\xi,1)\le v(\xi,x_1)$ for all $\xi<x_1$, which holds for $x_1\le x^*$, the root of (\ref{eq:xstar}) below, since then $v(\xi,x_1)\ge v_{\rm d}(x_1)\ge Q(x_1,1)\ge Q(\xi,1)$ by Lemma~\ref{lem:mon}(i).

\subsection{Beyond the rank-2 boundary}\label{sec:beyond}
For $m=2$ the value function is continuously differentiable in the open cone $\Omega=\{0<x_1<x_2\}$, across the boundaries $\beta_1,\beta_2$ included. The proof does not use smooth fit; it uses the ray representation (\ref{eq:ray}) and the fact that $x_2$ enters the right side of (\ref{eq:two}) only as the upper limit of integration.

\begin{proposition}[Regularity]\label{prop:reg}
In $\mathrm A_2$ and $\mathrm B_2$ the value function $v$ is continuously differentiable in $\Omega$, and $g=\partial_2v$ is differentiable along the rays and satisfies
\begin{equation}\label{eq:gode}
(x_1\partial_1+x_2\partial_2)\,g+(1+x_2)\,g=\min\{Q(x_2,2)-v,\ 0\}\qquad\text{in }\Omega .
\end{equation}
\end{proposition}
\begin{proof}
Let $P$ be the right side of (\ref{eq:two}). The integrands are bounded and continuous in $(\xi,\bx)$ except at $\xi=x_1$, a null set, so $P$ is continuous.

{\it Step 1.} For $\xi<x_2$ the state $\bx\cup\xi$ does not contain $x_2$, and the rank of $\xi$ does not depend on $x_2$; thus $x_2$ enters $P$ only as the upper limit. In $\Omega$ the integrand near $\xi=x_2$ is $\min\{Q(\xi,2),v(x_1,\xi)\}$, which is continuous in $\xi$, hence
$$\partial_2P(\bx)=\min\{Q(x_2,2),\,v(\bx)\},$$
a continuous function in $\Omega$.

{\it Step 2.} Differentiate (\ref{eq:ray}) in $x_2$ at fixed $x_1$. Since $\partial_{x_2}[P(r\bx)]=r\,\partial_2P(r\bx)$, and the $x_2$-derivative of the integrand is bounded and continuous in $r\in(0,1]$,
\begin{equation}\label{eq:dv}
\partial_2v(\bx)=-v(\bx)+e^{-x_2}\int_0^1e^{x_2r}\big[P(r\bx)+\partial_2P(r\bx)\big]\dd r ,
\end{equation}
so $\partial_2v$ exists and is continuous in $\Omega$.

{\it Step 3.} By (\ref{eq:ray}) the function $f(s)=v(s\bx)$ satisfies $(e^{sx_2}f)'=e^{sx_2}P(s\bx)/s$, whose right side is continuous for $s>0$. Hence $f$ is continuously differentiable, and
$$D v:=(x_1\partial_1+x_2\partial_2)v=P-x_2v$$
is continuous. In the coordinates $x_1=s$, $x_2=\kappa s$ the function $F(s,\kappa)=v(s,\kappa s)$ has the continuous partial derivatives $\partial_sF=(Dv)/s$ and $\partial_\kappa F=s\,\partial_2v$. So $F$ is continuously differentiable, and so is $v$, because $(x_1,x_2)\mapsto(s,\kappa)$ is a diffeomorphism on $\Omega$.

{\it Step 4.} Put $g=\partial_2v$ and substitute $r=u/s$ in (\ref{eq:dv}) at the point $s\bx$:
$$s\,e^{sx_2}\big(g+v\big)(s\bx)=\int_0^se^{x_2u}\big[P(u\bx)+\partial_2P(u\bx)\big]\dd u .$$
The right side is continuously differentiable in $s$. Differentiating, and using $\frac{\dd}{\dd s}v(s\bx)=(P-sx_2v)(s\bx)/s$ from Step~3, we obtain $s\,\frac{\dd}{\dd s}g(s\bx)+(1+sx_2)\,g(s\bx)+v(s\bx)=\partial_2P(s\bx)$. With Step~1 this is (\ref{eq:gode}) at the point $s\bx$.
\end{proof}

In the region $x_2>\beta_2(x_1)$ the $2$-records are not acceptable, so $Q(x_2,2)>v$ and the right side of (\ref{eq:gode}) vanishes: $(x_1\partial_1+x_2\partial_2)g+(1+x_2)g=0$. Along a ray $x_2=\kappa x_1$ this is $(x_2e^{x_2}g)'=0$, hence
\begin{equation}\label{eq:transport}
\frac{\partial v}{\partial x_2}(x_1,x_2)=C(\kappa)\,\frac{e^{-x_2}}{x_2},\qquad \kappa=\frac{x_2}{x_1},
\end{equation}
for a function $C$ of one variable. Since $\partial_2v$ is continuous (Proposition~\ref{prop:reg}), $C(\kappa)$ is the value of $x_2e^{x_2}\partial_2v$ computed from the closed form on $\mathcal R$ at the point where the ray meets the boundary. As $\kappa\to\infty$ this point tends to $(0,\beta_2(0))$ and $C(\infty)=-K_{\rm A}$ in $\mathrm A_2$, $-K_{\rm B}$ in $\mathrm B_2$, by Remark~\ref{rem:vone}. Integrating (\ref{eq:transport}) from the boundary, we obtain
\begin{equation}\label{eq:beyond}
v(x_1,x_2)=Q(\beta_2(x_1),2)+\int_{\beta_2(x_1)}^{x_2}C\Big(\frac\xi{x_1}\Big)\frac{e^{-\xi}}{\xi}\,\dd\xi,\qquad x_1<\gamma,\ x_2>\beta_2(x_1).
\end{equation}
For $x_1\ge\gamma$ no $2$-record is acceptable and the whole region $x_2>x_1$ is of this kind: $v(x_1,x_2)=v_{\rm d}(x_1)+\int_{x_1}^{x_2}C(\xi/x_1)e^{-\xi}\xi^{-1}\dd\xi$, where $v_{\rm d}(x)=v(x,x)$ is the value on the diagonal. On the diagonal the second integral in (\ref{eq:two}) vanishes, and
\begin{equation}\label{eq:diag}
v_{\rm d}'(x)+v_{\rm d}(x)=\frac1x\int_0^x\min\{Q(\xi,1),v(\xi,x)\}\,\dd\xi,\qquad v_{\rm d}(\gamma)=Q(\gamma,2),
\end{equation}
and $v(\binf)=\lim_{x\to\infty}v_{\rm d}(x)$.

\paragraph{The diagonal critical score.} As long as all records are acceptable the right side of (\ref{eq:diag}) is $\frac1x\int_0^xQ(\xi,1)\dd\xi$, and the equation is linear, with the same solution as on $\mathcal R$. Since $Q(\xi,1)-v(\xi,x)$ is increasing in $\xi$, the worst record has $\xi=x$, and all records are acceptable up to the root $x^*$ of the balance equation
\begin{equation}\label{eq:xstar}
Q(x,1)=v_{\rm d}(x).
\end{equation}
For $x_1>x^*$ the boundary $\beta_1$ bends away from the diagonal. Since $\beta_1(x_1)\le x^*$ for every $x_1$, the root of (\ref{eq:rho}) needs the value $v(\beta_1,x_1)=v_{\rm d}(\beta_1)+\int_{\beta_1}^{x_1}C(\xi/\beta_1)e^{-\xi}\xi^{-1}\dd\xi$ with $\beta_1<x^*$ only. Here (\ref{eq:transport}) applies on the whole segment from $(\beta_1,\beta_1)$ to $(\beta_1,x_1)$, because $\beta_1>\gamma$ by Lemma~\ref{lem:b1gamma}. 
 Therefore the {\it critical scores} $\beta_2$, $\beta_1$ and the initial threshold $H_0=\beta_1(\infty)$ are all given by explicit equations. The initial threshold is itself a balance point: a record of level $H_0$ arriving at time $0$ into the empty history leads to the state $(H_0,\infty)$, and $H_0$ is the root of $Q(H_0,1)=v(H_0,\infty)$, where $v(x,\infty):=\lim_{x_2\to\infty}v(x,x_2)$; this is (\ref{eq:rho}) with $x_1=\infty$. 

Beyond $x^*$ the function $v_{\rm d}$ is not elementary, since (\ref{eq:diag}) is an integro-differential equation of Volterra type. One differentiation, however, removes the integral and leaves two quadratures of functions that are already known.

\begin{proposition}[Quadratures for ${\bf a}_2$ and ${\bf b}_2$]\label{prop:value}
In $\mathrm A_2$ and in $\mathrm B_2$ put $\Theta(y):=\int_y^1C(1/s)\,\dd s$ for $0<y\le1$, and $K:=-x^*e^{x^*}v_{\rm d}'(x^*)$. For $x>x^*$
\begin{equation}\label{eq:ode}
\big(x\,e^{x}\,v_{\rm d}'(x)\big)'=\Theta\Big(\frac{\beta_1(x)}x\Big),
\end{equation}
and
\begin{equation}\label{eq:value}
v(\binf)=R(x^*)+\int_{x^*}^{\infty}\Theta\Big(\frac{\beta_1(s)}s\Big)J_1(s)\,\dd s,\qquad R(x^*)=Q(x^*,1)-K\,J_1(x^*),
\end{equation}
where $R(x^*)$ is the loss of the record rule $\mathrm R_b$ with $b=x^*$, and $K=\tfrac12x^{*2}e^{x^*}$ in $\mathrm A_2$, $K=2e^{x^*}-x^{*2}-2x^*-2$ in $\mathrm B_2$.
\end{proposition}
\begin{proof}
Multiply (\ref{eq:diag}) by $x$ and differentiate, which is allowed by Proposition~\ref{prop:reg}. For $x>x^*$ the integrand $\min\{Q(\xi,1),v(\xi,x)\}$ equals $Q(\xi,1)$ for $\xi\le\beta_1(x)$ and $v(\xi,x)$ above; it is continuous at $\xi=\beta_1(x)$ and equals $v_{\rm d}(x)$ at $\xi=x$. Hence
$$(x\,v_{\rm d}')'+(x\,v_{\rm d})'=v_{\rm d}(x)+\int_{\beta_1(x)}^{x}\partial_2v(\xi,x)\,\dd\xi .$$
For $x\ge x^*$ every state $(\xi,x)$ with $\xi<x$ lies above the rank-2 boundary. Indeed, by Lemma~\ref{lem:shape} and Proposition~\ref{prop:edge} we have $\beta_2(\xi)\le\beta_2(0)$ for $\xi<\gamma$, and there is no boundary for $\xi\ge\gamma$; and $\beta_2(0)<x^*$: in $\mathrm A_2$ the increasing function $(1+b)e^b$ equals $3$ at $b_{\rm A}$ and $5$ at $x^*$, and in $\mathrm B_2$ the increasing function $J_2$ equals $1$ at $b_{\rm B}$ and $x^*+\frac12>1$ at $x^*$ (here $x^*>\frac12$, because $J_2(x)-x-\frac12$ is increasing and negative at $x=\frac12$, where $J_2(\frac12)=0.5702$). So no $2$-record is acceptable in these states, and (\ref{eq:transport}) gives $\partial_2v(\xi,x)=C(x/\xi)e^{-x}/x$ for all $\xi<x$. With $\xi=sx$ the integral equals $e^{-x}\Theta(\beta_1(x)/x)$. This is $(xv_{\rm d}')'+xv_{\rm d}'=e^{-x}\Theta$, that is (\ref{eq:ode}). Integrating from $x^*$ gives 
$$xe^xv_{\rm d}'(x)=-K+\int_{x^*}^x\Theta(\beta_1(s)/s)\dd s.$$ 
On $[0,x^*]$ the function $v_{\rm d}$ is the function $G$ of the record rule, so $v_{\rm d}(x^*)=Q(x^*,1)$ and, by Proposition~\ref{prop:Rb}, $R(x^*)=G(x^*)+x^*e^{x^*}J_1(x^*)G'(x^*)=Q(x^*,1)-KJ_1(x^*)$. Finally $v(\binf)=v_{\rm d}(x^*)+\int_{x^*}^\infty v_{\rm d}'(x)\dd x$; inserting $v_{\rm d}'$, dividing by $xe^x$ and exchanging  the order of integration in the double integral (the function $\Theta$ is bounded) gives (\ref{eq:value}). The values of $K$ follow from $v_{\rm d}(x)=\frac12+\frac x2+\frac52e^{-x}$ in $\mathrm A_2$ with $(1+x^*)e^{x^*}=5$, and from $v_{\rm d}(x)=1-e^{-x}(2J_2(x)-x)$ in $\mathrm B_2$ with $J_2(x^*)=x^*+\frac12$.
\end{proof}
Since $C<0$, the optimal value is the loss of the record rule with $b=x^*$ plus a negative integral. The numbers are
\begin{center}\small
\begin{tabular}{llll}
\toprule
problem & $R(x^*)$ & integral in (\ref{eq:value}) & value\\
\midrule
$\mathrm A_2$ & $1.668057$ & $-0.013716$ & ${\bf a}_2=1.654341$\\
$\mathrm B_2$ & $0.214435$ & $-0.014975$ & ${\bf b}_2=0.199459$\\
\bottomrule
\end{tabular}
\end{center}
with $K=1.146846$ in $\mathrm A_2$ and $K=0.842425$ in $\mathrm B_2$. In $\mathrm A_2$ the function $C$ is used in its parametric form (Section~\ref{sec:A2}); in $\mathrm B_2$ it is tabulated, and $\Theta$ is computed from the table, so ${\bf b}_2$ is obtained to six digits. The weight $\beta_1(s)$ is itself a root of the balance equation with $v_{\rm d}$ and $C$ on $[0,x^*]$ only. The values for $x_1>x^*$ are not the analytic continuation of the formula valid for $x_1<x^*$: in the ray representation (\ref{eq:ray}) the bend of $\beta_1$ at $x^*$ enters through the ray point $r=x^*/x_1$.
\section{Problem $\mathrm A_2$}\label{sec:A2}
Here $Q(\xi,1)=1+\xi$, $Q(\xi,2)=2+\xi$, $Q_\infty=3$, and $q(i)=\min(i,3)$. Lemma~\ref{lem:lb} gives $v(x_1,x_2)\ge1+e^{-x_1}+e^{-x_2}$. A record of level $\xi<x_1$ is acceptable provided $1+\xi\le v(\xi,x_1)$, which is guaranteed if $1+\xi\le1+e^{-\xi}+e^{-x_1}$; this holds for all $\xi<x_1\le L_{\rm A}$, where $L_{\rm A}=0.8526\ldots$ is the root of $xe^x=2$ (the worst case being $\xi=x_1$).

\begin{theorem}\label{thm:A2}
On $\mathcal R\cap\{x_1\le L_{\rm A}\}$
\begin{equation}\label{eq:A2closed}
v(x_1,x_2)=\frac32+\frac{x_2}2-\frac{x_1}{x_2}+\Big(\frac32+\frac{x_1}{x_2}\Big)e^{-x_2}.
\end{equation}
The diagonal point is $\gamma=0.388925\ldots$, which is the root of $(3+\gamma)e^{\gamma}=5$, so $\mathcal R\subset\{x_1\le L_{\rm A}\}$. The rank-$2$ boundary is the curve $x_1=\varphi(x_2)$ for
\begin{equation}\label{eq:phiA}
\varphi(x_2)=\frac{x_2\,(\,3-(1+x_2)e^{x_2}\,)}{2\,(e^{x_2}-1)},\qquad \gamma\le x_2\le b_{\rm A},
\end{equation}
which starts at the point $(0,b_{\rm A})$ of the edge $x_1=0$, $b_{\rm A}=0.617642\ldots$, and ends at the diagonal point $(\gamma,\gamma)$; as $x_2$ decreases from $b_{\rm A}$ to $\gamma$, $\varphi(x_2)$ increases from $0$ to $\gamma$. Equivalently, the boundary point on the ray $x_2=\kappa x_1$ is the root $x_2$ of $\kappa=2(e^{x_2}-1)/\big(3-(1+x_2)e^{x_2}\big)$, which is a strictly increasing function of $x_2$ from $\kappa=1$ at $x_2=\gamma$ to $\kappa=\infty$ at $x_2=b_{\rm A}$.
\end{theorem}
\begin{proof}
On this set (\ref{eq:lin}) holds with the right-hand side 
$$\int_0^{x_1}(1+\xi)\dd\xi+\int_{x_1}^{x_2}(2+\xi)\dd\xi=2x_2-x_1+\frac12x_2^2.$$ For $f(s)=v(s,\kappa s)$ this is $f'+\kappa f=2\kappa-1+\frac12\kappa^2s$, $f(0)=3$, whence
$$f(s)=\frac32-\frac1\kappa+\frac{\kappa s}2+\Big(\frac32+\frac1\kappa\Big)e^{-\kappa s},$$
which is (\ref{eq:A2closed}) with $x_2=\kappa s$, $x_1/x_2=1/\kappa$. The balance equation (\ref{eq:mu}), $2+x_2=v(x_1,x_2)$, becomes $\frac{1+x_2}2+\frac{x_1}{x_2}=\big(\frac32+\frac{x_1}{x_2}\big)e^{-x_2}$, which is solved for $x_1$ by (\ref{eq:phiA}). At $x_2=b_{\rm A}$ we have $(1+x_2)e^{x_2}=3$ and $\varphi=0$. On the diagonal $v(x,x)=\frac12+\frac x2+\frac52e^{-x}$, and the difference $Q(x,2)-v(x,x)=\frac32+\frac x2-\frac52e^{-x}$ vanishes exactly when $(3+x)e^x=5$, at $\gamma_*=0.388925\ldots<L_{\rm A}$. At $x=L_{\rm A}$ the difference equals $0.86$. Let $\gamma$ be as in Lemma~\ref{lem:shape}. If $\gamma>L_{\rm A}$, then the diagonal states with $x\le L_{\rm A}$ lie in $\mathcal R$, where (\ref{eq:A2closed}) holds, and the difference would be positive at $L_{\rm A}$ and hence have its zero below $L_{\rm A}$, a contradiction. So $\gamma\le L_{\rm A}$, (\ref{eq:A2closed}) holds up to the diagonal point, and $\gamma=\gamma_*$ by uniqueness of the zero.
\end{proof}
The closed form agrees on the edge with $v(0,x)=1+G(x)$, as it should by Proposition~\ref{prop:edge}. 
For the continuation region (\ref{eq:A2closed}) gives the function $C$ in a parametric form: at the boundary point $(x_1,x_2)=(\varphi(x_2),x_2)$ we have $\kappa=x_2/\varphi$ and
$$C(\kappa)=x_2e^{x_2}\Big(\frac12+\frac{\varphi}{x_2^2}\Big)-x_2\Big(\frac32+\frac\varphi{x_2}+\frac{\varphi}{x_2^2}\Big),$$
which tends to $-K_{\rm A}$ as $x_2\to b_{\rm A}$.

\paragraph{The diagonal and the rank-1 boundary.} As long as all records are acceptable, (\ref{eq:diag}) is linear and its solution is the same as on $\mathcal R$, which is $v_{\rm d}(x)=\frac12+\frac x2+\frac52e^{-x}$. The balance equation (\ref{eq:xstar}), $1+x=v_{\rm d}(x)$, is $(1+x)e^x=5$, with the root $x^*=0.944456\ldots$ (this is also the hyperbola constant of the record rule with $Q_\infty=3$). The rank-1 boundary is the root of $1+\beta_1=v_{\rm d}(\beta_1)+\int_{\beta_1}^{x_2}C(\xi/\beta_1)e^{-\xi}\xi^{-1}\dd\xi$ with $x_2$ in place of $x_1$; for $x_2\to\infty$ it gives the initial threshold $H_0=\beta_1(\infty)=0.90169$, the root of the balance equation $1+H_0=v(H_0,\infty)$ of the arrival at time $0$. 

By Proposition~\ref{prop:value} the value is ${\bf a}_2=1.654341\ldots$; Section~\ref{sec:num} gives independent checks. For comparison, consider the record rule $\mathrm R_b$ of Section~\ref{sec:one}, which stops at the first {\it record} below the hyperbola $x=b/\bar t$ and never accepts a $2$-record. With $Q(\xi,1)=1+\xi$ and $Q_\infty=3$, its best constant $b$ is the root of the balance equation $Q(b,1)=G(b)$, which is $(1+b)e^b=5$, that is $b=0.944456\ldots$ ($=x^*$), and the expected loss is $1.668057$ by Proposition~\ref{prop:Rb}. The optimal rule in $\mathrm A_2$ improves on this only by $0.8\%$, which is the integral in (\ref{eq:value}) divided by ${\bf a}_2$.

\begin{remark}[$\mathrm A_m$ where every rank is kept]\label{rem:Am}
Let the state be $x_1<\dots<x_m$ in $\mathrm A_m$ and put $S=x_1+\dots+x_{m-1}$. On the set of states where an arrival of every initial rank $i\le m$ is acceptable, equation (\ref{eq:bellman}) is linear,
$$\sum_{i=1}^mx_i\,\partial_iv+x_m\,v=\sum_{j=1}^m\int_{x_{j-1}}^{x_j}(j+\xi)\,\dd\xi=\tfrac12x_m^2+m\,x_m-S,\qquad x_0:=0,$$
and its solution with $v(\bxero)=m+1$ is
$$v(\bx)=\frac{x_m+2m-1}{2}-\frac{S}{x_m}+e^{-x_m}\Big(\frac32+\frac{S}{x_m}\Big).$$
For $m=2$ this is (\ref{eq:A2closed}). The lower levels enter through their sum only. The balance for the rank $m$, $Q(x_m,m)=v(\bx)$, becomes $S=\varphi(x_m)$ with the function $\varphi$ of (\ref{eq:phiA}), the same for every $m$. This surface reaches the corner $x_1=\dots=x_m$ of the state space at the root of $(x+2m-1)e^x=2m+1$, which is $\gamma$ for $m=2$, and $0.2817$, $0.2203$ for $m=3,4$, respectively. With all levels equal to $x$ the formula is the function $G(x)=\frac{x+1}2+\frac{2m+1}2e^{-x}$ of the record rule. The set where every rank is acceptable contains a neighbourhood of the origin, where $v$ is close to $m+1$ and the loss of stopping is at most $m+\xi$; for $m\ge3$ we have not determined how far it extends, and the formula is claimed only there. For $m=3$ the method of Proposition~\ref{prop:value} does not give ${\bf a}_3$ by quadratures, because the unknown is a function of two variables (the value on the states with $x_2=x_3$); numerically ${\bf a}_3=1.7677$, against $1.80407$ for the record rule $\mathrm R_b$, $b=1.1708$ the root of $(1+b)e^b=7$.
\end{remark}

\section{Problem $\mathrm B_2$}\label{sec:B2}
The problem $\mathrm B_2$ has a very similar structure.
Here $Q(\xi,1)=1-e^{-\xi}(1+\xi)$, $Q(\xi,2)=1-e^{-\xi}$, $Q_\infty=1$, and $q(i)=\mathbf 1(i>2)$. Lemma~\ref{lem:lb} gives $v(x_1,x_2)\ge e^{-x_2}$. A record of level $\xi<x_1$ is acceptable as soon as $Q(\xi,1)\le e^{-x_1}$, which holds for all $\xi<x_1\le L_{\rm B}$, where $L_{\rm B}=1.1462\ldots$ is the root of $(2+L)e^{-L}=1$.

\begin{theorem}\label{thm:B2}
On $\mathcal R\cap\{x_1\le L_{\rm B}\}$ the value is $v=1-\Psi$, where
\begin{equation}\label{eq:psi}
\Psi(x_1,x_2)=e^{-x_2}\Big[2J_2(x_2)-J_2(x_2-x_1)-\frac{x_1}{x_2-x_1}\big(e^{\,x_2-x_1}-1\big)\Big].
\end{equation}
The diagonal point is $\gamma=0.705926\ldots$, the root of $2J_2(\gamma)-\gamma=1$, and $\mathcal R\subset\{x_1\le L_{\rm B}\}$. The rank-$2$ boundary $x_2=\beta_2(x_1)$ is given by the balance equation $\Psi(x_1,x_2)=e^{-x_2}$, that is, with $w=x_2-x_1$,
\begin{equation}\label{eq:muB}
2J_2(x_2)-J_2(w)-\frac{x_2-w}{w}\,(e^{w}-1)=1,
\end{equation}
which gives $\beta_2(0)=b_{\rm B}=0.804352\ldots$ (at $w=x_2$, where it reduces to $J_2(x_2)=1$) and the diagonal point $x_2=\gamma$ at $w\to0$. On the ray $x_2=\kappa x_1$ the boundary point is the root $x_2$ of $$2J_2(x_2)-J_2\big(x_2\tfrac{\kappa-1}\kappa\big)-\frac{1}{\kappa-1}\big(e^{x_2(\kappa-1)/\kappa}-1\big)=1.$$
\end{theorem}
\begin{proof}
On this set the right-hand side of (\ref{eq:lin}) is
$$\int_0^{x_1}\big(1-e^{-\xi}(1+\xi)\big)\dd\xi+\int_{x_1}^{x_2}\big(1-e^{-\xi}\big)\dd\xi=x_2-2+(1+x_1)e^{-x_1}+e^{-x_2}.$$
For $f(s)=v(s,\kappa s)$, equation (\ref{eq:lin}) after multiplication by $e^{\kappa s}/s$ is
$$(e^{\kappa s}f)'=e^{\kappa s}\Big(\kappa+\frac{-2+(1+s)e^{-s}+e^{-\kappa s}}{s}\Big),\qquad e^{\kappa s}f\to1\ \ (s\to0).$$
The right-hand side is integrable at $0$; integrating and expressing the result through $J_2$ gives (\ref{eq:psi}) with $x_2=\kappa s$, $x_1/x_2=1/\kappa$. (The substitution of (\ref{eq:psi}) into (\ref{eq:lin}) was also checked numerically, the residual being $10^{-11}$, and at $x_1=0$ it reduces to $1-e^{-x_2}J_2(x_2)$ of (\ref{GB}).) Letting $x_2\downarrow x_1=x$ gives $v(x,x)=1-e^{-x}(2J_2(x)-x)$, so $Q(x,2)-v(x,x)=e^{-x}(2J_2(x)-x-1)$, whose zero is $\gamma_*=0.7059\ldots<L_{\rm B}$, and which is positive at $L_{\rm B}$. As in the proof of Theorem~\ref{thm:A2}, Lemma~\ref{lem:shape} gives $\gamma=\gamma_*$.
\end{proof}
The value is ${\bf b}_2=0.19946\ldots$ (Proposition~\ref{prop:value}), thus the maximal probability of choosing one of the two best is about $1-{\bf b}_2=0.80054$. The function $C$ is obtained from (\ref{eq:psi}) as for $\mathrm A_2$, with $C(\infty)=-K_{\rm B}$.

\begin{remark}[The shape of $\beta_2$ in $\mathrm B_2$]\label{rem:beta2B}
Put $g(w)=(e^w-1)/w$, so that $J_2'=g$. Expanding (\ref{eq:muB}) in $x_1$ at fixed $x_2$, the terms of first order cancel and the equation becomes $J_2(x_2)=1-\frac12x_1^2g'(x_2)+O(x_1^3)$. Hence $\beta_2(x_1)=b_{\rm B}-\kappa x_1^2+O(x_1^3)$ with $\kappa=g'(b_{\rm B})/(2g(b_{\rm B}))=0.2832$: the boundary leaves the edge with zero slope. At the diagonal point the slope is $\dd x_2/\dd x_1=-0.2273$. In $\mathrm A_2$ the curve (\ref{eq:phiA}) leaves the edge with the slope $-0.5700$ and meets the diagonal with the slope $-0.6065$.
\end{remark}

\paragraph{The diagonal and the rank-1 boundary.} As long as all records are acceptable, $v_{\rm d}(x)=1-e^{-x}(2J_2(x)-x)$, and the balance equation (\ref{eq:xstar}) becomes $J_2(x)=x+\frac12$, with the root $x^*=1.218477\ldots$. The quadrature for the rank-1 boundary gives the initial threshold  $H_0=\beta_1(\infty)=1.0883\ldots$, the root of the balance equation $Q(H_0,1)=v(H_0,\infty)$ of the arrival at time $0$.

\begin{figure}[!t]
\centering
\input{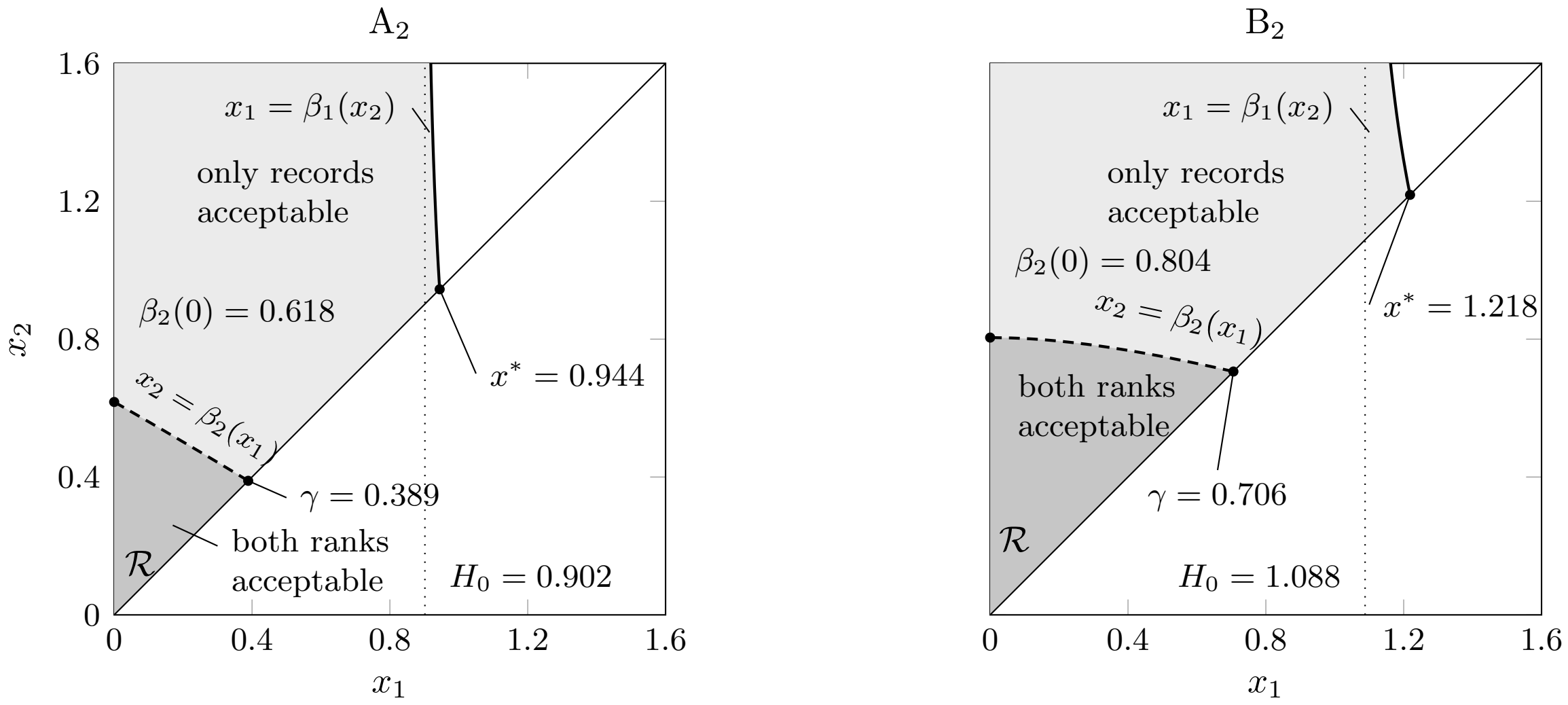}
\caption{Decision regions in the state plane $(x_1,x_2)$, $x_1<x_2$, in box coordinates, for $\mathrm A_2$ (left) and $\mathrm B_2$ (right). A record of level $\xi$ arriving in the state $(x_1,x_2)$ leads to the new state $(\xi,x_1)$, and a $2$-record to $(x_1,\xi)$; the shade of a state says what the optimal strategy does with the atom which leads to it. Dark grey, region $\mathcal R$: both ranks are acceptable; its upper boundary (dashed) is the curve $x_2=\beta_2(x_1)$ of balance $Q(x_2,2)=v(x_1,x_2)$, joining $(0,\beta_2(0))$ to $(\gamma,\gamma)$. Light grey: only records are acceptable, up to the curve $x_1=\beta_1(x_2)$ (solid) of balance $Q(x_1,1)=v(x_1,x_2)$, which coincides with the diagonal up to $(x^*,x^*)$ and then continues as a strip narrowing to $\beta_1(\infty)=H_0$ (dotted). White: continuation. The curve $\beta_2$ is the explicit curve (\ref{eq:phiA}) in $\mathrm A_2$ and the solution set of (\ref{eq:muB}) in $\mathrm B_2$; $\beta_1$ is computed by quadrature.}\label{fig:stop}
\end{figure}

\subsection{Stopping probabilities}
Let $p_i(\bx)$ be the probability that the optimal rule started in $\bx$ stops at a $i$-record. It solves the linear version of (\ref{eq:two}) for the optimal policy of Figure~\ref{fig:stop}: at an arrival the minimum is replaced by $\mathbf 1\{i=\text{rank}\}$ if the arrival is acceptable, and by $p_i$ at the new state otherwise, with the edge conditions $p_1(0,x)=0$ and $p_2(0,x)=1-e^{-x}$ for $x\le b_{\rm B}$, $p_2(0,x)=1-e^{-b_{\rm B}}+b_{\rm B}J_1(b_{\rm B})-b_{\rm B}J_1(x)$ for $x\ge b_{\rm B}$.
\begin{center}
\begin{tabular}{lccc}
\toprule
 & stops at a record & stops at a $2$-record & no stop ($\tau=1$)\\
\midrule
solution of the linear equation & $0.865$ & $0.041$ & $0.094$\\
simulation, $2.8\times10^7$ trajectories & $0.8646\pm0.0001$ & $0.0415\pm0.0001$ & $0.0939\pm0.0001$\\
\bottomrule
\end{tabular}
\end{center}
The probability of success is $1-{\bf b}_2=0.8005$, of which $0.7729$ comes from records and $0.0276\pm0.0001$ from $2$-records (simulation). All numbers in this subsection are numerical; the grid solution of the linear equation converges at first order only, because the boundary is located to within one grid step, so its entries are accurate to about $10^{-3}$. Thus $2$-records account for only $4.6\%$ of the stops. The probability $0.094$ is that of $\tau=1$: no suitable arrival occurred.

\subsection{A sample path of the threshold}
\begin{figure}[ht]
\centering
\input{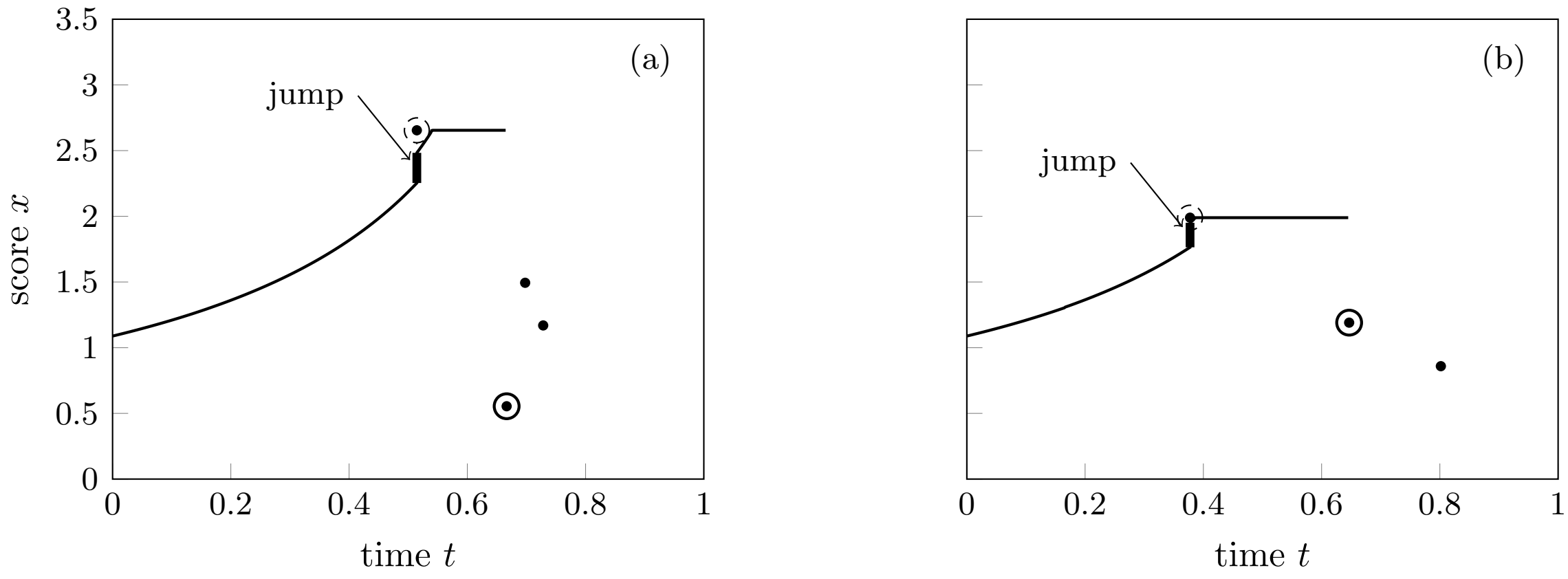}
\caption{Two realisations of $\Pi$ (unit intensity; only the atoms with $x\le3.5$ are shown) with the adaptive threshold $H_t$ (solid) of the optimal rule in $\mathrm B_2$. The first atom below the curve is selected (circled). Between arrivals the threshold drifts up. When a record arrives above the curve (dashed circle), the threshold jumps up (heavy segment) and then rises to the score of this record, where it stays flat. The realisations were chosen to show a jump; the median stopping time among the realisations with a stop is $0.42$.}\label{fig:path}
\end{figure}

Figure~\ref{fig:path} shows two realisations of $\Pi$ with the adaptive threshold $(H_t)$ of the optimal rule, defined as in Section~\ref{sec:mon} on the whole interval, with the initial value $H_0=1.09$. It drifts up; when a record arrives above the curve, it jumps up; after a further drift it reaches the lowest score $y^\circ_1$ and stays flat there (Proposition~\ref{prop:mon}). The rule stops at the first atom below the curve. The $2$-record boundary matters only if $\bar ty^\circ_1<\gamma=0.706$.

\section{Simple strategies}\label{sec:thr}
The optimal rules need tracking the $m$ lowest scores to compute
an adapted threshold.
We compare them with much simpler rules. A strategy which stops at the first atom below a nonrandom curve $x=h(t)$ is called a {\it memoryless} rule \cite{BF, Memoryless}. 
 In this section we derive the formulas; the numerical values are collected in Section~\ref{sec:num}. For $\mathrm B_m$ we evaluate the {\it failure} probability, the quantity which is minimised;
 $\pi_m(z)=\P[\text{Poisson}(z)\le m-1]$ denotes the probability of at most $m-1$ atoms in a region of area $z$.

In the rank-information version of choosing one of the $m$ best \cite{FS, GM, GZ} the failure probability of the optimal rule decays exponentially fast in $m$ (see Section~4.1 of \cite{Samuels}), while a rule with a single cutoff, accepting the first atom with initial rank at most $m$ arriving after time $c=m^{-1/(m-1)}$, has the success probability $c$ for $m\ge2$ \cite{GnedinStrategies} (eq.~(16), p.~162).
 For the full-information problem the simple rules below make the failure probability much smaller: e.g.\ for $m=2$ it is between $0.20$ and $0.27$, against $0.4264$ for the optimal rule with rank information \cite{FS} and $1-c=0.5$ for the single cutoff; the comparison for $m\le5$ is in Section~\ref{sec:num}. The success probability $c$ of the single cutoff is the maximum over the cutoff $t$ of the limiting success $r(t-t^r)/(r-1)$ given in (3.2) of \cite{FS}, with $r=m$.

\subsection{The record rule}
The {\it record rule} $\mathrm R_b$ of Section~\ref{sec:one} stops at the first {\it record} below the hyperbola $x=b/\bar t$; the atoms of the initial rank above $1$ are never accepted, whether or not they are below the hyperbola. For $\mathrm R_b$ with the loss $\mathrm B_m$ the computation of Section~\ref{sec:one} applies verbatim, with $g(\xi)=Q(\xi,1)=1-\pi_m(\xi)$ and $Q_\infty=1$: the loss of a record depends on its level only. By Proposition~\ref{prop:Rb} the failure probability is $G(b)+be^bJ_1(b)G'(b)$, with $G$ defined by (\ref{eq:G}). The same applies to $\mathrm A_m$, with $g(\xi)=1+\xi$ and $Q_\infty=m+1$; the optimal constant $b$ solves the balance equation $Q(b,1)=G(b)$, which for $\mathrm A_2$ is $(1+b)e^b=5$.

\subsection{Memoryless rules: a shifted hyperbola}
Let $h(t)=c+b/\bar t$ with $b>0$, $c+b>0$, and let $\tau$ be the first atom below the curve $x=h(t)$, whatever its rank. This is a memoryless rule. The success probability in $\mathrm B_m$ is
$$
W_m(c,b)=\int_0^\infty\dd x\int_{t_x}^1 e^{-ct}\bar t^{\,b}\,\pi_m\big(\Lambda(x)+x\bar t\big)\,\dd t,
$$
where $t_x=0$, $\Lambda(x)=0$ for $x\le c+b$, and for $x>c+b$
$$
t_x=1-\frac{b}{x-c},\qquad \Lambda(x)=x-c-b-b\log\frac{x-c}{b}.
$$
The failure probability is $1-W_m(c,b)$. Indeed, the first atom below the curve has the density $e^{-ct}\bar t^{\,b}$ in $\dd x\,\dd t$ for $x<h(t)$, that is for $t>t_x$, since $\int_0^th(s)\,\dd s=ct-b\log\bar t$. The atoms with score below $x$ which arrived before $t$ are above the curve, in the lens $\{(s,y):s<t,\ h(s)\le y<x\}$ of area $\Lambda(x)$, and the atoms after $t$ below $x$ number Poisson$(x\bar t)$, all independent. The final rank is at most $m$ if and only if the total number is at most $m-1$. The case $c=0$ is the first atom below $b/\bar t$.

\subsection{Memoryless rules: a constant threshold}\label{sec:flat}
For $\mathrm B_m$ consider also the constant curve $h(t)\equiv b$: the rule stops at the first atom with score below $b$. Since the strip $(0,t)\times(0,b)$ is empty at the time $t$ of stopping, this atom is a record, and no smaller score has arrived before it. The pair $(t,x)$ has the density $e^{-bt}\dd t\,\dd x$ on $[0,1]\times[0,b]$, the box coordinate is $\xi=x\bar t$, and with probability $e^{-b}$ no atom arrives and the loss is $Q_\infty=1$. Thus the failure probability equals $G(b)$ of (\ref{eq:G}) with $g(\xi)=Q(\xi,1)=1-\pi_m(\xi)$: this is the same as starting the greedy phase from the state $\xi=b$. In closed form,
\begin{equation}\label{flatB}
1+e^{-b}\Big(\sum_{j=1}^{m-1}\frac{(m-j)b^j}{j\cdot j!}-m\,J_2(b)\Big),
\end{equation}
where we used the identity $\int_0^a\big(1-\pi_m(z)\big)\dd z=\E\,(N_a-m)_+=a-m+e^{-a}\sum_{j<m}(m-j)a^j/j!$ for $N_a$ Poisson$(a)$, and then integrated in $t$ with $s=\bar t$. For $m=1$ the expression is $1-e^{-b}J_2(b)$, which is (\ref{GB}). The numerical values are in Section~\ref{sec:num}.

\subsection{Rank-wise hyperbolas}
For $m\ge2$ let the rule accept an $i$-record, $i\le m$, if and only if its box coordinate is below $b_i$, that is, below the hyperbola $b_i/\bar t$. The rule is stationary in box coordinates, and its loss is computed from (\ref{eq:bellman}) with the minimum replaced by the policy. For $m=2$ this is a two-dimensional problem, and for $m=3$ a three-dimensional one (Section~\ref{sec:num}).

\section{Numerical results and comparison}\label{sec:num}
Numerical values are rounded to the last digit shown, and simulated values are given with a standard error.

\paragraph{Computed from the equations.} The values ${\bf a}_2,{\bf b}_2$ are given by the quadratures of Proposition~\ref{prop:value}: ${\bf a}_2=1.654341$ and ${\bf b}_2=0.199459$ (the latter from the tabulated function $C$, to six digits). They were also found by marching the optimality equation (\ref{eq:two}) in the top component $x_2$, on a grid $x_j=x_{\min}e^{kh}$ in both coordinates, so that the shrinking $\bx\mapsto(1-\dd t)\bx$ is a diagonal shift of the grid; the right side involves $v$ only at states with a smaller top component, so no iteration over the free boundaries is needed. With the steps $h=0.01,0.005,0.0025$ the observed order is $1.98$, and the extrapolation (Aitken) gives ${\bf a}_2=1.65435$ and ${\bf b}_2=0.19946$, within $2\times10^{-5}$ of the quadratures; a march of the diagonal equation (\ref{eq:diag}) alone gives $1.654341$ for $\mathrm A_2$. The scheme reproduces the closed forms (\ref{eq:A2closed}) and (\ref{eq:psi}) to $8\times10^{-5}$ and $3\times10^{-6}$, the rank-2 boundary (\ref{eq:phiA}) to $10^{-4}$, and the rank-1 boundary, given by the quadrature of Section~\ref{sec:two}, to $3\times10^{-5}$ ($\mathrm A_2$) and $10^{-6}$ ($\mathrm B_2$); the quadrature also gives $H_0=0.90169$ and $1.0883$. The losses of the rules of Section~\ref{sec:thr} are values of the integrals there, optimised over the parameters by numerical quadrature and minimisation; all the numerics were done in Python (NumPy, SciPy), and the columns $b/\bar t$ and constant $b$ were recomputed with independent code, which agrees to $10^{-6}$; the formula (\ref{flatB}) was also checked by direct numerical integration and, for $m=2$, by simulation. For the rank-wise hyperbolas the march is applied to the policy equation; for $m=3$ it is three-dimensional, and the steps $0.04,0.02,0.01$ with extrapolation give $0.09837$ (simulation of the rule, $8\times10^6$ paths: $0.0984\pm0.0001$); this is the computed value of one rule, known to about $10^{-4}$, so it bounds ${\bf b}_3$ only up to this error.

\paragraph{Computed by simulation.} The stopping probabilities of Section~\ref{sec:B2} ($2.8\times10^7$ trajectories of $\Pi$ with the optimal rule, with the boundaries $\beta_2$ from (\ref{eq:muB}) and $\beta_1$ from the quadrature), and the sample paths of Figure~\ref{fig:path}. In $7\times10^6$ simulated realisations the median stopping time of the optimal rule in $\mathrm B_2$ among the realisations with a stop is $0.422$, with the quartiles $0.209$ and $0.655$; if no stop is counted as $\tau=1$ the three numbers are $0.467$, $0.230$ and $0.739$. No atom is accepted ($\tau=1$) with the frequency $0.0939\pm0.0001$, in agreement with $0.094$ of the linear equation.

\paragraph{Simple strategies.} The table gives the failure probabilities for $\mathrm B_m$; the last column is the optimal value ${\bf b}_m$, which for $m=3$ is a numerical estimate, see below.
\begin{center}\footnotesize
\setlength{\tabcolsep}{2.5pt}
\begin{tabular}{lllllll}
\toprule
$m$ & memoryless & memoryless & memoryless & record rule & rank-wise & optimal\\
& $b/\bar t$ & $c+b/\bar t$ & constant $b$ & $b/\bar t$ & hyperbolas & ${\bf b}_m$\\
\midrule
$1$ & $0.462363$ & $0.444463$ & $0.482649$ & $0.419836$ & -- & $0.419836$\\
 & $b=0.7415$ & $c=0.8188,\,b=0.2903$ & $b=1.503$ & $b=0.8044$ & & \\[1mm]
$2$ & $0.238279$ & $0.222758$ & $0.273530$ & $0.214435$ & $0.19971$ & $0.19946$\\
 & $b=1.0544$ & $c=0.9680,\,b=0.4838$ & $b=2.018$ & $b=1.2185$ & $b_1=1.127,\,b_2=0.744$ & \\[1mm]
$3$ & $0.128339$ & $0.116775$ & $0.164514$ & $0.118438$ & $0.098$ & $\approx0.098$\\
 & $b=1.3624$ & $c=1.115,\,b=0.672$ & $b=2.493$ & $b=1.6223$ & $b_{1,2,3}=1.402,\,1.135,\,0.716$ & \\
\bottomrule
\end{tabular}
\end{center}

For $m=3$ the optimum of the rule $c+b/\bar t$ over $(c,b)$ is flat along a ridge. The excess of the best of the single-curve rules of the table over the optimal failure probability is $0\%$, $7.5\%$ and about $20\%$ for $m=1,2,3$ (the record rule for $m=1,2$ and the rule $c+b/\bar t$ for $m=3$; the record rule for $m=1$ is optimal), and the rank-wise hyperbolas reduce the excess at $m=3$ to less than $1\%$. For $m=2$ the optimal rule differs from the record rule by the region of $2$-records of Figure~\ref{fig:stop} and by the rank-$1$ boundary beyond $x^*$, which lies below the constant level $x^*$ of the record rule; the effect is $0.015$; two constants $b_1,b_2$ recover all of it but $2.5\times10^{-4}$.

\paragraph{Rank information against full information.} The following table gives the success probabilities for choosing one of the $m$ best. The entries for rank information are limiting values: the single cutoff $c=m^{-1/(m-1)}$ (the success probability is $c$, see Section~\ref{sec:thr}) and the optimal rule of \cite{FS} (Table~1, success $=1-Q_m$). The entries for full information are one minus the failure probabilities of the memoryless rules (constant $b$ and $b/\bar t$ of this section, with the optimal $b$ in brackets), and the optimal values of \cite{GMal} for $n=2500$, which we take as the limits.
\begin{center}\small
\begin{tabular}{lccccc}
\toprule
 & \multicolumn{2}{c}{rank information} & \multicolumn{3}{c}{full information}\\
\cmidrule(lr){2-3}\cmidrule(lr){4-6}
$m$ & single cutoff & optimal \cite{FS} & constant $b$ & $b/\bar t$ & optimal \cite{GMal}\\
\midrule
$2$ & $0.5000$ & $0.5736$ & $0.7265$ $(2.018)$ & $0.7617$ $(1.054)$ & $0.8007$\\
$3$ & $0.5774$ & $0.7082$ & $0.8355$ $(2.493)$ & $0.8717$ $(1.362)$ & $0.9025$\\
$4$ & $0.6300$ & $0.7987$ & $0.8979$ $(2.947)$ & $0.9291$ $(1.670)$ & $0.9513$\\
$5$ & $0.6687$ & $0.8603$ & $0.9353$ $(3.386)$ & $0.9601$ $(1.979)$ & $0.9750$\\
\bottomrule
\end{tabular}
\end{center}
The full-information optimal entries for $m=4,5$ are less accurate than for $m\le3$: the grids in \cite{GMal} have $56$ and $32$ intervals, and the column for $m=5$ is not monotone in $n$ at $n=750$, so we regard them as reliable to about $10^{-3}$ only. A change of $10^{-3}$ in the optimal failure probability moves the two excesses given below for $m=5$ over the ranges $149$--$170\%$ and $53$--$66\%$, and those for $m=4$ over $105$--$114\%$ and $43$--$49\%$.

The single level $b$ is the simplest rule, it uses neither the time nor the history of the scores, and yet its success probability exceeds that of the optimal rule with rank information by $0.153$, $0.127$, $0.099$ and $0.075$ for $m=2,\dots,5$ (the gain shrinks since the rank-information optimum approaches $1$ exponentially fast). Compared with the full-information optimum its success probability is lower by $0.074$, $0.067$, $0.053$ and $0.040$, which in the failure probability is an excess of $37\%$, $69\%$, $110\%$ and $159\%$, whereas the rule $b/\bar t$ has the excess of $20\%$, $32\%$, $46\%$ and $60\%$ (for $m=2,\dots,5$, relative to the entries of \cite{GMal}; relative to ${\bf b}_2$ and ${\bf b}_3$ of the text the excesses for $m=2,3$ differ in the first decimal only). For $m=1$ the excess of the single level is $15\%$. Thus in the examples considered the constant level is attractive for its simplicity, but its relative excess grows with $m$. The rules which use the ranks only through a few constants capture most of the value: the record rule is within $7.5\%$ for $m=2$ and the rank-wise hyperbolas within $1\%$ for $m=2,3$. For $\mathrm A_2$ the record rule $\mathrm R_b$ has the loss $1.668057$, against ${\bf a}_2=1.65434$, within $0.8\%$.

\paragraph{Comparison with the literature.} The table compares our values with those of \cite{BF} for the truncations of Robbins' problem ($\mathrm A_m$), and with those of \cite{GMal} for choosing one of the $m$ best ($\mathrm B_m$).
\begin{center}\small
\begin{tabular}{lllll}
\toprule
 & ours & \cite{BF} & \cite{GMal} & comment\\
\midrule
${\bf a}_1$ & $1.462488$ \,(analytic) & $1.462$ & -- & agree\\
${\bf a}_2$ & $1.65434$ \,(quadrature) & $1.658$ & -- & supported by recursion\\
${\bf a}_3$ & $1.7677$ \,(numerical) & $1.782$ & -- & see below\\
${\bf b}_2$ & $0.19946$ \,(numerical) & -- & $0.1993$ at $n=2500$ & extrapolation agrees\\
${\bf b}_3$ & $\approx0.098$ \,(numerical) & -- & $0.0975$ at $n=2500$ & see below\\
${\bf b}_4$ & -- & -- & $0.0487$ at $n=2500$ & not computed here\\
${\bf b}_5$ & -- & -- & $0.0250$ at $n=2500$ & not computed here\\
\bottomrule
\end{tabular}
\end{center}

{Goldenshluger and Malinovsky \cite{GMal}:} The optimal $n$-step success probabilities for the problem of choosing one of the $m$ best, computed by the recursion of \cite{GMal} (Table~1 there), are for $m=2$: $0.8044$, $0.8009$, $0.8007$ at $n=100,1000,2500$, for $m=3$: $0.9026$, $0.9025$ at $n=1000,2500$, for $m=4$: $0.9517$, $0.9513$ and for $m=5$: $0.9758$, $0.9750$ (four digits). They decrease in $n$ like $1/n$ (for $m=1$ the entries $0.5804$ and $0.5803$ at $n=1000,2500$ against the limit $0.58016$). The entries we use are those at $n=2500$, which differ from the limits by $O(1/n)$. The extrapolation gives $0.8005$ for $m=2$, in agreement with $1-{\bf b}_2=0.80054$, and $0.9024\pm0.0002$ for $m=3$, that is ${\bf b}_3=0.0976\pm0.0002$ (the entry $0.9025$ at $n=2500$ is slightly above the limit). An independent dynamic program of ours for $\mathrm B_{n,3}$ (the three lowest scores on a grid, extrapolated in the grid size; checked with three grid shapes) gives the values $0.09056$, $0.09422$, $0.09605$ at $n=50,100,200$, whose increments halve, hence ${\bf b}_3=0.0979\pm0.0002$. The two estimates differ by $3\times10^{-4}$; we state ${\bf b}_3\approx0.098$, and the percentages for $m=3$ above are accurate to one unit. 
For $m=4,5$ the grids of \cite{GMal} are coarse (see above), and we use the entries at $n=2500$ without extrapolation.

{Bruss and Ferguson \cite{BF}:} The values $W^{(m)}$ of \cite{BF} (Section~4, Table~3) are extrapolated from the optimal $n$-step values for $n\le17$; ${\bf a}_1$ is confirmed to the three digits given. For $\mathrm A_2$ the discrepancy is $0.0037$ ($0.2\%$). Since the finite-$n$ values increase with $n$ and are bounded by ${\bf a}_2$, a value above ${\bf a}_2$ can only come from an overshoot of the extrapolation. To check this independently we computed the exact $n$-step dynamic program for the loss $\min(i,3)+(n-k)X_k$, with the state the two lowest scores, on grids of $1500$ and $3000$ points (the program reproduces the  values $V_2=1.25$ and $V_3=1.3915$ of the rank loss). The values for $n=50,100,200$ are $1.638947$, $1.646645$, $1.650495$ (Richardson-extrapolated in the grid size), the differences decrease by the factor $2$, and the Richardson extrapolations $2W_{2n}-W_n$ are $1.654343$ and $1.654345$, in agreement with ${\bf a}_2=1.65434$ and not with $1.658$. For $m=3$ the value $1.782$ of \cite{BF} is $0.014$ above our ${\bf a}_3=1.7677$, which is supported by three computations of different kinds (agreeing to $10^{-6}$) and by an independent finite-$n$ program ($1.7675\pm0.0005$).


\begin{thebibliography}{99}
\bibitem{Bruss} F.\,T.\ Bruss, What is known about Robbins' problem? {\it J.\ Appl.\ Probab.} {\bf 42} (2005), 108--120.
\bibitem{BF} F.\,T.\ Bruss and T.\,S.\ Ferguson, Minimizing the expected rank with full information, \emph{J.\ Appl.\ Probab.} {\bf 30} (1993), 616--626.
\bibitem{BFhp} F.\,T.\ Bruss and T.\,S.\ Ferguson, Half-prophets and Robbins' problem of minimizing the expected rank, in: \emph{Athens Conference on Applied Probability and Time Series}, Lecture Notes in Statistics {\bf 114}, Springer, 1996, 1--17.
\bibitem{CMRS} Y.\,S.\ Chow, S.\ Moriguti, H.\ Robbins and S.\,M.\ Samuels, Optimal selection based on relative rank (the ``secretary problem''), \emph{Israel J.\ Math.} {\bf 2} (1964), 81--90.
\bibitem{RSC} Y.\,S.\ Chow, H.\ Robbins and D.\ Siegmund, \emph{Great Expectations: The Theory of Optimal Stopping}, Houghton Mifflin, 1971.
\bibitem{Swan} R.\ Dendievel and Y.\ Swan, One step further: an explicit solution to Robbins' problem when $n=4$, \emph{Mathematica Applicanda} {\bf 44} (2016), 135--148.
\bibitem{DY} E.\,B.\ Dynkin and A.\,A.\ Yushkevich, \emph{Theorems and Problems on Markov Processes}, Plenum Press, 1969.
\bibitem{FS} A.\,Q.\ Frank and S.\,M.\ Samuels, On an optimal stopping problem of Gusein-Zade, \emph{Stoch.\ Process.\ Appl.} {\bf 10} (1980), 299--311.
\bibitem{GiSamuels} J.\ Gianini and S.\,M.\ Samuels, The infinite secretary problem, \emph{Ann.\ Probab.} {\bf 4} (1976), 418--432.
\bibitem{GM} J.\ Gilbert and F.\ Mosteller, Recognizing the maximum of a sequence, \emph{J.\ Amer.\ Statist.\ Assoc.} {\bf 61} (1966), 35--73.
\bibitem{GnedinStrategies} A.\ Gnedin, Multicriteria extensions of the best choice problem: sequential selection without linear order, in: \emph{Strategies for Sequential Search and Selection in Real Time}, {\it Contemporary Mathematics} {\bf 125} (1992), 152--172.
\bibitem{GFI} A.\ Gnedin, On the full-information best-choice problem, \emph{J.\ Appl.\ Probab.} {\bf 33} (1996), 678--687.
\bibitem{BC-SPA} A.\ Gnedin, Best choice from the planar Poisson process, \emph{Stoch.\ Process.\ Appl.} {\bf 111} (2004), 317--354.
\bibitem{GnedinJAP} A.\ Gnedin, Optimal stopping with rank-dependent loss, \emph{J.\ Appl.\ Probab.} {\bf 44} (2007), 996--1011.
\bibitem{GnedinCorners} A.\ Gnedin, Corners and records of the Poisson process in quadrant, \emph{Electron.\ Commun.\ Probab.} {\bf 13} (2008), 187--193.
\bibitem{GnedinUniform} A.\ Gnedin, Optimal stopping for the uniform distribution,  arXiv:2601.04371 (2026).
\bibitem{Memoryless} A.\ Gnedin and M.\,C.\,S.\ Carreira, The memoryless best-choice problem, arXiv:2607.16145 (2026).
\bibitem{GKS} A.\ Gnedin, P.\ Kozie{\l} and M.\ Sulkowska, Running minimum in the best-choice problem, \emph{Extremes} {\bf 26} (2023), 157--182.
\bibitem{GMal} A.\ Goldenshluger and Y.\ Malinovsky, On optimal stopping problems with full information and rank-dependent reward, preprint (2026).
\bibitem{GZ} S.\,M.\ Gusein-Zade, The problem of choice and the optimal stopping rule for a sequence of independent trials, \emph{Theory Probab.\ Appl.} {\bf 11} (1966), 472--476.
\bibitem{Mucci} A.\,G.\ Mucci, Differential equations and optimal choice problems, \emph{Ann.\ Statist.} {\bf 1} (1973), 104--113.
\bibitem{PS} G.\ Peskir and A.\ Shiryaev, \emph{Optimal Stopping and Free-Boundary Problems}, Birkh\"auser, 2006.
\bibitem{Rusch} L.\ R\"uschendorf, Approximative solutions of optimal stopping and selection problems, \emph{Mathematica Applicanda} {\bf 44} (2016), 17--44.
\bibitem{Samuels} S.\,M.\ Samuels, Secretary problems, in \emph{Handbook of Sequential Analysis}, Dekker, 1991, 381--405.
\bibitem{Shiryaev} A.\,N.\ Shiryaev, \emph{Optimal Stopping Rules}, Springer, 1978.
\end{thebibliography}
\end{document}